\documentclass[10pt]{article}
 \usepackage[dvips]{graphicx}
 \usepackage{amsfonts, amssymb, amsmath, amsthm}
\newcommand \RR {\mathbb{R}}
 \newcommand \eps \varepsilon

\begin{document}

\title{\bf A version of the Glimm method\\
based on generalized Riemann problems}

\date{May 17, 2006} 
\maketitle
\begin{center}

\author{\large John M. Hong\footnote{Department of Mathematics,
National Central University, Chung-Li 32054, Taiwan.
\\
\noindent
E-mail: jhong@math.ncu.edu.tw}
\, and \, Philippe G. LeFloch\footnote{Laboratoire J.-L. Lions \& C.N.R.S.,
University of Paris VI, 75252 Paris, France.
\\
\noindent E-mail : lefloch@ann.jussieu.fr
 }
}
\vspace{0.3cm}

\end{center}
\medskip

 \begin{abstract}
 \noindent
We introduce a generalization of Glimm's random choice method, which provides us with an
approximation of entropy solutions to quasilinear hyperbolic system of balance laws.
The flux-function and the source term of the equations may depend
on the unknown as well as on the time and space variables. The method is
based on local approximate solutions of the generalized Riemann problem, which form building blocks in our scheme
and allow us to take into account naturally the effects of the flux and source terms. 
To establish the nonlinear stability of these approximations, we investigate nonlinear interactions between generalized
wave patterns.
This analysis leads us to a global existence result for quasilinear hyperbolic systems with source-term, and 
applies, for instance, to the compressible Euler equations in general geometries and to
hyperbolic systems posed on a Lorentzian manifold.
\end{abstract}

\

 %
 %
 \setcounter{section}{1}
 \noindent{\bf \large{1. Introduction}} \\
 \setcounter{equation}{0}

 \noindent

{\bf 1.1 Hyperbolic systems of balance laws.} This paper\footnote{This paper is based on notes written
by the second author in July 1990, 
for his Habilitation memoir (Chapter~XI) at the University of Paris VI.}
 is concerned with the approximation of entropy solutions to the
Cauchy problem for a quasilinear hyperbolic system
 \begin{eqnarray}
   \partial_t u + \partial_x f (t,x,u) = g (t, x, u), \quad t > 0, \; x \in \RR,
   \label{1.1}
 \end{eqnarray}
 \begin{eqnarray}
 u (0, x) = u_0 (x), \quad x \in \RR,
 \label{1.2}
 \end{eqnarray}
where $u=u(t,x) \in \RR^p$ is the unknown.
We propose here a generalized version of the Glimm scheme \cite{Gl} 
which allows us to deal with a large class of mappings $f,g$ and take into account the geometric effect 
of the flux and source terms. 
Our scheme is based on an approximate solver for the {\sl generalized Riemann problem}, 
based on an asymptotic expansion introduced by LeFloch and Raviart \cite{Le-Ra}.
The approach provides high accuracy and stability, under mild restrictions on the equation and the data. 

In \eqref{1.1}, the flux $f=f(t,x,u) \in \RR^p$ and the
source-term $g=g(t,x,u) \in \RR^p$ are given smooth maps defined
for all $(t,x,u) \in \RR_+
\times
\RR
\times {\cal U}$, where ${\cal U}$ is a small neighborhood of the origin in $\RR^p$, and the initial data $u_0 :
\RR
\to {\cal U}$ is a function with bounded total variation.
We assume that the Jacobian
matrix $A (t, x, u) := \frac {Df}{Du} (t, x, u)$ admits $p$
real and distinct eigenvalues,
$$
\lambda_1 (t, x, u) < \lambda_2 (t, x, u) < \ldots < \lambda_p (t, x, u),
$$
and therefore a basis of right-eigenvectors $r_j (t,x,u)$ ($1 \leq j \leq p$),
Finally, we assume that each characteristic field is either genuinely nonlinear
($\nabla \lambda_j (t, x, u) \cdot  r_i (t, x, u) \ne 0$)
or linearly degenerate
($\nabla \lambda_j (t, x, u) \cdot r_i (t, x, u) = 0$).

One important motivation for considering general balance laws \eqref{1.1} comes from the theory of general relativity.
In this context, the vector $u$ typically consists of fluid variables as well as (first order derivatives)
of the coefficients of an unknown, Lorentzian metric tensor. (See \cite{BLSS, BLF} and the reference therein.)
One can also freeze the metric coefficients and concentrate on the dynamics of the fluid.  
For instance, the compressible Euler equations describing the dynamics of a gas flow in general geometry
read:
\begin{eqnarray}
\begin{split}
\label{1.10}
&\partial_t \rho+\partial_x(\rho v)=-\frac {\partial_x a}{a} \rho
v- \frac{\partial_t a}{a}\rho,\\
&\partial_t (\rho v)+\partial_x(\rho v^2+p)=-\frac {\partial_x
a}{a}
(\rho v^2) - \frac{\partial_t a}{a} \, \rho v,\\
&\partial_t (\rho E)+\partial_x(\rho v E+pv)=-\frac {\partial_x
a}{a}
(\rho v E +p v) - \frac{\partial_t a}{a} \, \rho E\\
\end{split}
\end{eqnarray}
where $a=a(t,x)>0$ can be regarded as the cross section of a time-dependent (moving) duct, and $\rho$, $v$,
$p(\rho,e)$, $e$, and $E=e + u^2/2$ are the density, velocity, pressure, internal energy, and total energy of the gas, respectively.
The system \eqref{1.10} describes a situation where the fluid does not affect the variation of the
duct; i.e.~the function $a(t,x)$ is given and, for simplicity, smooth. 
The system \eqref{1.10} is of the form \eqref{1.1} with $ u=(\rho,
\rho v, \rho E)^T, \ f=f(u)=(\rho v, \rho v^2+p,
\rho v E+pv)^T$ and $g=g(t,x,u)=- \frac{\partial_x a}{a}g_1(u)- \frac{\partial_t a}{a}u$ where
$g_1(u) = (\rho v, \rho v^2,
\rho v E+pv)^T$. 

We are interested in solutions to \eqref{1.1}-\eqref{1.2} which have bounded total variation in space for all times
and satisfy the equations in the sense of distributions, together with an entropy condition \cite{La2, Dafermos,LeFloch-book}.
In the special case that
$$
f = f (u),   \quad g = 0,
$$
the existence of global entropy solutions was established by Glimm \cite{Gl}, 
assuming that the initial data $u_0(x)$ has sufficiently small total variation.
Recall that two main ingredients in Glimm's random choice method are
(1) the solutions of Riemann problems and
(2) a projection step based on a sequence of randomly chosen points.

Let us first indicate some of the earlier work on the subject. The system \eqref{1.1} with
$$
f = f (x,u),  \quad g = g(x, u),
$$
was treated in pioneering work by Liu \cite{Li2,Li4}, via a suitable extension of the Glimm method:
the approximate solutions are defined by pasting together steady state solutions, i.e.,
solutions $v=v(x)$ of the ordinary differential equation
 \begin{eqnarray*}
\frac {d}{dx} \big( f (x,v) \big)= g (x, v).
\label{1.5}
\end{eqnarray*}
He established the existence of solutions defined in a finite
interval of time $[0,T)$ as long as either $T$ or the $L^1$ norms
of $g$ and $\partial g / \partial u$ are sufficiently small.
Next, assuming in addition that the eigenvalues of the matrix
$A(x,u)$ never vanish (so that no resonance takes place), Liu
deduced a global existence result (with $T=+\infty$). Steady-state solutions were also used in the 
work by Glimm, Marshall, and Plohr \cite{GMP}.

For more general mappings $f,g$, the existence for \eqref{1.1}-\eqref{1.2} is established by
Dafermos and Hsiao \cite{Da-Hs} and Dafermos \cite{Dafermos, Dafermos2}. 
They assume that $f_x(u^*,t,x) = g(u^*,t,x)=0$
at some (equilibrium) constant state $u^*$, hence $u^*$ is a solution of \eqref{1.1}
around which \eqref{1.1} can be formally linearized. They also require that
the linearized system satisfies a dissipative property.
Their main result concerns the consistency and stability of a generalization of 
the Glimm method, yielding therefore the 
global existence of entropy solutions to \eqref{1.1}. In
\cite{Da-Hs}, the approximate solutions to the Cauchy problem on each time
step are based on classical Riemann solutions with initial data suitably modified by both the source term $g$
and the map $\theta := A^{-1} \, f_x $.

Next, Amadori et al. \cite{AmadoriGuerra, Am-Go-Gu} developed further techniques to
establish the existence of solutions for a large class of systems having $f=f(u)$ and $g=g(x,u)$, 
and discussed Dafermos-Hsiao dissipative condition. 
For some particular systems (of two or three equations) the condition that the total variation
be small can be relaxed; see for instance Luskin and Temple \cite{Lu-Te}, Groah and Temple \cite{GT},
Barnes, LeFloch, Schmidt, and Stewart \cite{BLSS}, and the references cited therein.
In these papers,  the decreasing of a total variation functional (measured with
respect to a suitable chosen coordinate) was the key to establish the stability of the scheme.

\

{\bf 1.2 A new version of the Glimm method.}
In the present paper
we provide an alternative approach to Dafermos-Hsiao's method,
and introduce a generalized version of the Glimm scheme for general mappings $f,g$.
Integrability assumptions will be required (and discussed later on) 
on the matrix $A$ and the mapping $q : \RR_+ \times \RR \times {\cal U} \to \RR^p$
defined by
\begin{eqnarray}
q (t, x, u) := g (t, x, u) - \frac{\partial f}{\partial x} (t,x, u).
\label{1.6}
\end{eqnarray} 
It should be emphasized that only this combination of the source and the flux will be important in our approach, 
which can be summarized as follows.

First, we study the generalized Riemann problem associated
with the system \eqref{1.1}, i.e.~the Cauchy problem with
piecewise constant initial data. The existence of solutions
defined locally in spacetime in a neighborhood of the initial
discontinuity was studied in Li and Yu \cite{Li-Yu} and
Harabetian \cite{Ha}. Contrary to the case where $f,g$ only
depend upon the unknown $u$, no closed formula is available for
the solutions of the generalized Riemann problem. We propose here
an {\sl approximate Riemann solver}, inspired by a technique of asymptotic
expansion introduced by Ben-Artzi and Falcovitz \cite{BM} (for the gas dynamics equations)
and LeFloch and Raviart \cite{Le-Ra} (for general hyperbolic systems of balance laws); see also \cite{Bo-Le-Ra}.

Our scheme for solving approximately the generalized Riemann
problem can be re-interpreted as a splitting algorithm (the
hyperbolic operator and the source term being decoupled). Since
an approximate (rather than an exact) solution to the generalized
Riemann problem is used, it is crucial to establish an error
estimate which we achieve in Proposition~2.1 below, under a mild
assumption on the data $u_0,f,g$. This estimate will be necessary to ensure the consistency of our generalized Glimm
method.

Second, 
we study the nonlinear interaction of waves issuing from two generalized Riemann problems, 
and establish a suitable extension of Glimm's estimates \cite{Gl} to the general system \eqref{1.1}; cf.~Proposition 3.3.
This is a key, technical part of our analysis.  

Third, we introduce our scheme and prove its stability in total variation, under the assumption that
the initial data $u_0$ has sufficiently small total variation and that the total amplification 
due to (the derivatives of) $f,g$ to the total variation of the solution is sufficiently small; cf.~Theorem 4.3.
More precisely, we impose that  
$$
\frac {\partial^2 A}{\partial t \partial u}, \quad
\frac {\partial^2 A}{\partial x \partial u},
\quad
q,
\quad
\frac{\partial q}{\partial u}  
$$
are sufficiently small in $L^1 (\RR_+ \times \RR)$.

Finally, we conclude with the convergence of the proposed scheme (Cf.~Theorem 5.1) 
which yields the global existence of entropy solutions for the Cauchy problem
\eqref{1.1}-\eqref{1.2}. The solution satisfies an entropy inequality and has bounded total variation
in $x$ for all $t \geq 0$. Our results cover in particular the case
 \begin{eqnarray}
\partial_t u + \partial_x f (u) = g (t), 
\label{1.7}
 \end{eqnarray}
for which global existence of entropy solutions is established under the sole assumption
 \begin{eqnarray}
\int^{+\infty}_{0} |g(t)| \, d t <<1.
\label{1.8}
 \end{eqnarray}
Without further restriction on the flux $f$, 
this condition is clearly necessary in order to apply the Glimm method, since, 
for instance in the trivial case $p = 1$ and $f = 0$, \eqref{1.7} reduces to the
differential equation
 \begin{eqnarray}
\partial_t u = g (t).
\label{1.9}
 \end{eqnarray}

On one hand, the condition \eqref{1.8} holds if and only if every solution of
\eqref{1.9} remains close to a constant state, which is a necessary
condition in order to apply the Glimm method. On the other
hand, when one of the eigenvalues of the system \eqref{1.1} vanishes, the
amplitude of solutions could become arbitrarily large and the
solutions would not remain bounded ---except when
the source term satisfies a ``damping" property in time.

As a direct application, the global existence of entropy solutions to \eqref{1.10} follows,
if the source $g$ and its derivative $\frac{\partial g}{\partial u}$
are sufficiently small in $L^1(\RR_+ \times \RR)$, which is the case, for instance, 
if the support of $(a_t, a_x)$ is sufficiently small.
 
 \


 \setcounter{section}{2}
 \noindent{ \bf \large {2. An approximate solver for the generalized Riemann problem}} \\
 \setcounter{equation}{0}

In the present section, we introduce an approximate solution to the
generalized Riemann problem associated with the system \eqref{1.1}, and we
derive an error estimates (see Proposition~2.1 below).

Given $t_0> 0$, $x_0\in \RR$, and two constant states $u_L, u_R \in \RR^p$, we consider the {\sl generalized Riemann
problem}, denoted by $R_G (u_L, u_R ; t_0, x_0)$, and consisting of the following equations and initial conditions:
\begin{eqnarray}
\partial_t u + \partial_x f (t, x, u) = g (t, x, u), \quad t > t_0, \; x \in \RR,
\label{2.1}
\end{eqnarray}
 \begin{eqnarray}
u (0, x) =\left\{
\begin{array}{ll} u_L, \quad x < x_0,\\
 u_R, \quad x > x_0.
\end{array}
\right.
\label{2.2}
 \end{eqnarray}

Replacing $f$ and $g$ in \eqref{2.1} by $f (t_0, x_0, u)$ and
$0$, respectively, the problem $R_G (u_L, u_R ; t_0, x_0)$
reduces to the {\sl classical Riemann problem},which we denote by
$R_C (u_L, u_R ;t_0, x_0)$, that is the equations
\begin{eqnarray}
\partial_t u + \partial_x f (t_0, x_0, u) = 0, \quad u (t, x)
\in \RR^p, \quad t > t_0, \quad x \in \RR
\label{2.3}
\end{eqnarray}
together with the initial data (2.2). This problem was solved by Lax under the assumption that the initial jump
$|u_R - u_L|$ be sufficiently small:
the solution to $R_C (u_L, u_R ;t_0, x_0)$ is self-similar
(i.e.~depends only on $\frac {x - x_0}{t - t_0}$) and consists of
at most $(p + 1)$ constant states $u_L = u_0, u_1,\ldots, u_p =
u_R$, separated by rarefaction waves, shock waves or contact
discontinuities; see Figure 2.1.

The following terminology and notation will be used throughout this paper.
Let $W_C = W_C (\xi ; u_L, u_R ; t_0, x_0)$ be the solution of $R_C (u_L, u_R ;t_0, x_0)$ with
$\xi = (x - x_0)/(t - t_0)$. We say that the problem $R_C (u_L, u_R ; t_0, x_0)$ is solved
by the elementary waves $(u_{i- 1}, u_i)$ ($i = 1,\ldots,p$)
if each constant state $u_i$ belongs to the {\sl i-wave curve} ${\cal
W}_i (u_{i - 1})$ issued from the state $u_{i - 1}$ in the phase
space, and $(u_{i- 1}, u_i)$ is called an $i-wave$ of $R_C (u_L,
u_R ; t_0, x_0)$. When the $i$-characteristic field is
genuinely nonlinear, the curve ${\cal W}_i (u_{i -1})$ consists
of two parts, the $i$-rarefaction curve and the $i$-shock curve
issuing from $u_{i -1}$;
if $i$-characteristic field is
linearly degenerate, the curve ${\cal W}_i (u_{i - 1})$ is a
$C^2$ curve of $i$-contact discontinuities. Call
$\eps_i$ the strength of the $i$-wave $(u_{i- 1}, u_i)$ along the $i$-curve,
so that, for a genuinely nonlinear $i$-field, we can assume that
$\eps_i \geq 0$  if  $(u_{i-1},u_i)$ is a rarefaction wave, and
$\eps_i \leq 0$   if   $(u_{i-1},u_i)$ is a shock wave. On the other hand,
$\eps_i$ has no specific sign if $(u_{i-1},u_i)$ is a contact discontinuity.

Let $\eps_i (u_L, u_R ; t_0, x_0)$ denote the {\sl wave
strength} of the $i$-wave $(u_{i - 1}, u_i)$ in the Riemann
problem $R_C (u_L, u_R ; t_0, x_0)$, and vector $\eps = (\eps_1,
\ldots,\eps_p)$ denote the {\sl wave strength} of $R_C
(u_L, u_R ; t_0, x_0)$ (so $|\eps |$ is equivalent to the
total variation of $W_C (\xi ; u_L, u_R ; t_0, x_0)$). In
addition, we let $\sigma^-_i
=
\lambda_i (u_{i - 1}, t_0, x_0)$ and $\sigma^+_i =  \lambda_i
(u_i, t_0, x_0)$ be the lower and upper speeds of the
$i$-rarefaction wave $(u_{i - 1}, u_i)$ respectively, and
$\sigma_i$ be the speed of the $i$-shock or $i$-contact
discontinuity. If the $i$-wave is a shock or a contact
discontinuity we set $\sigma^-_i =
\sigma^+_i = \sigma_i$.\\

From the implicit function theorem we deduce that the states
$u_i$ and the speeds $\sigma^\pm_i$ are smooth functions of $u_L$,
 $u_R$, $t_0$, and $x_0$. Moreover, one can check that
$u_i = u_L + O(1) |u_R - u_L|$ ($i = 0,1,\ldots,p$),
and, for an $i$-shock $(u_{i - 1}, u_i)$,
$$
\sigma_i = \lambda_i (u_{i - 1}; t_0, x_0) + O(1) \, |u_i - u_{i - 1}|, \quad i = 1,2,\ldots,p,
$$
where $O(1)$ is bounded function possibly depending on $u_L, u_R \in {\cal U}$, $t_0\geq 0$, and $x_0 \in \RR$.\\

Consider next the generalized Riemann problem on which a large
literature is available \cite{Li-Yu, Ha, BM, Bo-Le-Ra,Le-Ra}. First, we recall \cite{Li-Yu} that the
solution of $R_G (u_L, u_R ; t_0, x_0)$ is piecewise smooth
and has a local structure which is similar to the one of the associated classical Riemann problem $R_C (u_L, u_R ; t_0, x_0)$.
Following \cite{Le-Ra} we consider an {\sl approximate Riemann solution}
of the problem $R_G (u_L, u_R ; t_0, x_0)$,  denoted by $W_G (t, x ;
u_L, u_R ; t_0, x_0)$ and  defined by
\begin{eqnarray}
W_G (t, x ; u_L, u_R ; t_0, x_0) = W_C (\xi) + (t - t_0) \, q(t_0, x_0, W_C (\xi))
\label{2.7}
\end{eqnarray}
for $t > t_0$ and $x \in \RR$. Here, the function $q(t,x,u)$ is given by
\eqref{1.6}, and
$$
\xi = \frac {x - x_0}{t -t_0}, \quad
W_C (\xi) = W_C (\xi ; u_L, u_R ; t_0, x_0). $$
Observe that the function
$W_G (t, x ; u_L, u_R ; t_0, x_0)$ is constructed as a
superposition of the corresponding classical Riemann solution $W_C
(\xi ; u_L, u_R ; t_0, x_0)$ and an asymptotic expansion term $(t
- t_0) q (t_0, x_0, W_C (\xi))$ (see Figure 2.2).

Within a region where function $W_C(\xi)$ is a constant, the
function $W_G(t,x ; u_L, u_R ; t_0, x_0)$ is a
linear function of $t$, namely,
\begin{eqnarray}
W_G (t, x; u_L, u_R ; t_0, x_0) = u_i + (t - t_0) q (t_0, x_0,
u_i),
\quad  \sigma^+_i < \frac {x}{t} < \sigma^-_{i + 1}
\label{2.8}
\end{eqnarray}
for $i = 0, 1,\ldots, p$. By convention, $\sigma^+_0 := - \infty$ and
$\sigma^-_{p + 1} := + \infty$. Whenever there will be no ambiguity, we will
use the notation $W_G(t,x)$ or $W_G (t, x ; u_L, u_R)$ for $W_G
(t,x ; u_L, u_R ; t_0, x_0)$.

To describe the structure of $W_G
(t, x ; u_L, u_R ; t_0, x_0)$, it is convenient to say that the approximate solution
$W_G (t, x, u_L, u_R ; t_0, x_0)$ consists of an $i$-wave $(u_{i
- 1}, u_i)$ if $(u_{i - 1}, u_i)$ is an $i$-wave of the
corresponding classical Riemann solution $W_C (\xi ; u_L, u_R ; t_0, x_0)$.\\

\begin{center}
\includegraphics[height = 6cm, width=10 cm]{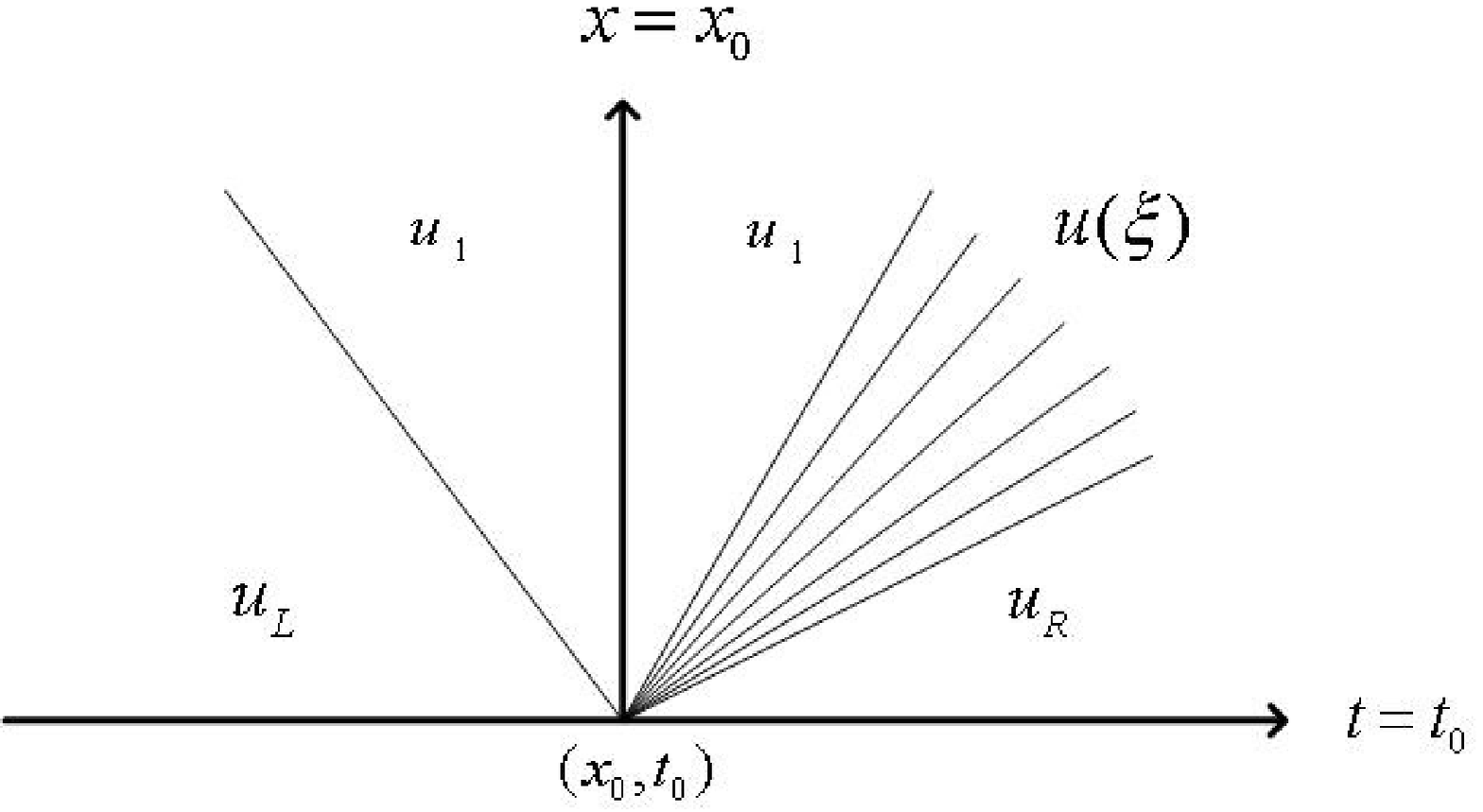} \\
  Figure 2.1 :  Classical Riemann solution (p=2), \\
  $u_L$, $u_1$ and $u_R$ are constant states, $u=u(\xi)$ is a
  function of $\xi =\frac{x}{t}$.
\end{center}
\begin{center}
\includegraphics[height = 6cm, width=10 cm]{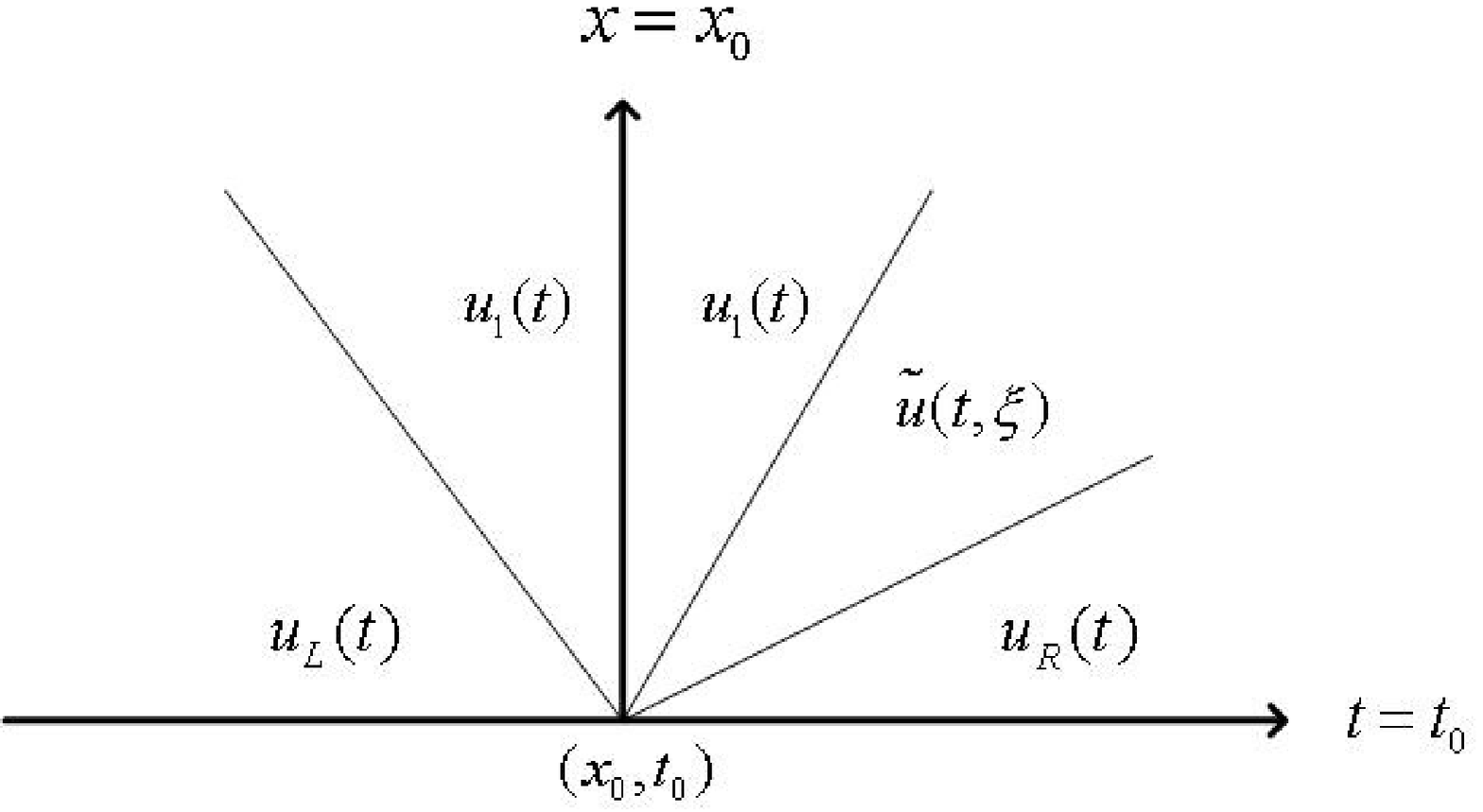} \\
  Figure 2.2 :  Generalized Riemann solution (p=2),\\
  $u_{L}(t)$,
  $u_{1}(t)$,
  $u_{R}(t)$ are functions of $t$ and
  $\tilde{u}(t,\xi)$ is constructed by \eqref{2.7}.
\end{center}

We now prove that the function $W_G(t,x)$ defined in \eqref{2.7}
approximately solves the problem $R_G (u_L, u_R ; t_0, x_0)$,
by evaluating the discrepancy between $W_G(t,x)$ and the exact solution of $R_G (u_L,
u_R ; t_0, x_0)$. Given any $s > 0$ and $r > 0$, and any $C^1$ function $\theta :
\RR_+ \times \RR
\to \RR$ with compact support, we now show that the term
\begin{eqnarray}
\Delta(s, r ; \theta) := \int^{t_0 + s}_{t_0}
\int^{x_0 + r}_{x_0 - r}
\{W_G \, \partial_t \theta + f (t, x, W_G) \, \partial_x \theta +
g (t, x, W_G) \, \theta\Big) \, d x d t
\label{2.9}
\end{eqnarray}
is of third order in $r,s$, provided that the condition
\eqref{2.10} holds.

\vskip 0,5cm

\noindent
{\bf Proposition 2.1.} {\em \quad Let $\theta : \RR_+ \times \RR \to \RR$
be a compactly supported, $C^1$ function.
Then, for every $ (t_0, x_0) \in \RR_+ \times \RR$,
$u_L, u_R \in {\cal U} $, and any  positive numbers $s,r$
satisfying the (Courant-Friedrichs-Levy -type) stability condition
\begin{eqnarray}
\frac {s}{r} \sup |\lambda_i (t, x, u)| \leq 1
\label{2.10}
\end{eqnarray}
(the supremum being taken over $ \ 1 \leq i \leq p$, $ (t, x)
\in \RR_+ \times \RR$, and  u
$\in {\cal U} $), the \\function $
\ W_G (t, x) = W_G (t, x ; u_L, u_R ; t_0, x_0) \ $
satisfies
\begin{eqnarray}
\begin{split}
\Delta (s, r ; \theta)\; =& \;
  \int^{x_0 + r}_{x_0 - r}W_G (t_0 + s, \cdot) \theta (t_0 + s, \cdot) d x -\int^{x_0 + r}_{x_0 - r} W_G (t_0, \cdot) \theta(t_0,\cdot) d x\\
&+\int^{t_0 + s}_{t_0} f (\cdot, x_0 + r, W_G (\cdot, x_0 + r))\theta (\cdot, x_0 + r) d t \\
&-\int^{t_0 + s}_{t_0}f (\cdot, x_0 - r, W_G (\cdot, x_0 - r)) \theta(\cdot, x_0 - r) d t \\
&+ O(1) (s^2 + r^2) (s + r + |u_R - u_L|) ||\theta||_{{\cal C}^1},
\end{split}
\label{2.11}
\end{eqnarray}
\noindent
where $\Delta (s, r ; \theta)$ is given in \eqref{2.9} and
$
||\theta||_{{\cal C}^1} = ||\theta||_ {{\cal C}^0}
+ ||\partial_t \theta||_ {{\cal C}^0} + ||\partial_x \theta||_{{\cal C}^0}$.
}

The left-hand side of \eqref{2.11} vanishes when
$W_G(t,x)$ is replaced by the exact solution of $R_G(u_L, u_R ; t_0, x_0)$. Thus, the right hand side of
\eqref{2.11} represents the error due to the choice of
approximate solution $W_G(t,x)$.

\

\noindent{\bf Remark 2.2.}
1. Condition \eqref{2.10} ensures
that the waves in $R_C (u_L, u_R ; t_0, x_0)$ can not reach the
lines $\big\{x = x_0 \pm r \big\}$ for $t \leq t_0 + s$, so that the waves
in the rectangle region $D_{(t_o,x_0)} \equiv [x_0-r, x_0+r]\times
[t_0, t_0+s)$ do not interact with the waves outside $D_{(t_o,x_0)}$.

2. In a different context, Liu \cite{Li2} derived earlier an estimate similar to \eqref{2.10},
but for an approximation based on steady state solutions of the hyperbolic system and
with initial data consisting of two steady state solutions of (2.1) (with $f = f (u)$ and $g = g (x, u))$.

3. Our formula (2.11) yields a possible generalization to
the class of quasilinear systems (1.1) of the notion of (classical) Riemann
solver introduced by Harten and Lax in \cite{Ha-La}.

4. One can check similarly that $W_G$ satisfies an entropy inequality associated with an entropy pair
(when available). The error terms are completely similar to those found in (2.11).
This will be used to show that the weak solution generated by the random 
choice method satisfies all the entropy inequalities.

\vskip 0,5cm

\begin{proof}
\quad Without loss of generality, we can assume that $(t_0, x_0) = (0, 0)$.
Given a $C^1$ function $\theta $ with compact support in $\RR_+
\times \RR$, we define $m (t, x):= W_G \partial_t
\theta + f (t, x, W_G)
\partial_x
\theta + g (t, x, W_G) \theta$. From \eqref{2.9} we
have
 $\Delta (s, r ; \theta) =
\int^s_0 \int^r_{- r} m (t, x) d x d t$. Next, we decompose $\Delta
(s, r ; \theta)$ as
\begin{eqnarray}
\Delta (s, r ; \theta) ={\sum^p_{i = 0}}\Delta^1_i (s, r ; \theta) +\sum_{i-rare. \atop waves}\Delta^2_i (s, r ; \theta)
\label{2.12}
\end{eqnarray}
where
 \begin{eqnarray*}
\Delta^1_i (s, r ; \theta) := \int^s_0
\int^{\sigma^-_{i + 1} t}_{\sigma^+_i t}  m (t, x) \, dxdt, \quad
1 \leq i \leq p - 1,
 \end{eqnarray*}
 \begin{eqnarray*}
\Delta^1_0 (s, r ; \theta) := \int^s_0
\int^{\sigma^-_{1} t}_{- r}  m (t, x) \, dxdt,
\qquad \quad
\Delta^1_p (s, r ; \theta) := \int^s_0
\int^r_{\sigma^+_{p} t}  m (t, x) \, dxdt,
 \end{eqnarray*}
and (if the $i$-wave, $1 \leq i \leq p$, is a rarefaction wave)
 \begin{eqnarray*}
\Delta^2_i (s, r ; \theta) := \int^s_0
\int^{\sigma^+_i t}_{\sigma^-_i t} m(t, x) \, dxdt
\end{eqnarray*}\\

We first compute $\Delta^1_i$ in the region where classical
Riemann solution $W_C$ is a constant state. According to the form
of $W_G(t,x)$ in
\eqref{2.8}, it follows that
\begin{eqnarray}
W_G (t, x) = u_i + t \, q (0, 0 ; u_i)
\label{2.13}
\end{eqnarray}
for $\frac{x}{t} \in [\sigma^+_i, \sigma^-_{i + 1}],$ $i$ $\in$
$\big\{1, 2, \ldots, p-1 \big\}$. By a simple calculation and the definition of $q$ in \eqref{1.6}, we have
$$
\partial_t W_G + \partial_x f (t, x, W_G) - g (t, x, W_G) =
q (0, 0 ; u_i) - q (t, x, W_G)
$$
for $i$ $\in$ $\big\{1, 2,
\ldots, p-1 \big\}$.
By multiplication by the function $\theta$ and then using integration by parts, we obtain
\begin{eqnarray}
\begin{split}
\Delta^1_i (s, r ; \theta) \;=\; &\int^{\sigma^-_{i + 1} s}_{\sigma^+_i s}  W_G (s, x) \theta (s, x) d x  \\
&+ \int^s_0 (f (t, \sigma^-_{i + 1} t, W_G (t, \sigma^-_{i + 1} t))- \sigma^-_{i + 1} W_G (t, \sigma^-_{i + 1} t))  \theta (t, \sigma^-_{i +1} t) d t \\
&- \int^s_0 (f (t, \sigma^+_i t, W_G (t, \sigma^+_i t)) -\sigma^+_i W_G(t, \sigma^+_i t)) \theta (t, \sigma^+_i t) d t \\
&- \int^s_0 \int^{\sigma^-_{i + 1} t}_{\sigma^+_i t} \big( q (0, 0; u_i) - q (t, x, W_G) \big) \, \theta (t, x) \, dx d t.
\end{split}
\label{2.15}
\end{eqnarray}
\noindent
By the property that $q$ is Lipschitz continuous with respect to
$t$, $x$ and $u$ on the compact set $[0, s] \times [-r, r]$ and
the form of $W_G(t,x)$ in \eqref{2.13}, the last term on the right
hand side of
\eqref{2.15} can be estimated by $O(s^3) \,  ||\theta||_{{\cal
C}^0}$ with the bound $O(1)$ depending on $q$. Therefore,
equality (2.11) leads to
\begin{eqnarray}
\begin{split}
\Delta^1_i (s, r ; \theta) \;&=\; \int^{\sigma^-_{i + 1} s}_{\sigma^+_i s}W_G (s, x) \theta (s, x) d x\\
&+\displaystyle\int^s_0 (f (t, \sigma^-_{i + 1} t,W_G(t,\sigma^-_{i+1} t)) -\sigma^-_{i + 1}W_G (t, \sigma^-_{i + 1} t  - )) \theta (t, \sigma^-_{i + 1} t ) d t\\
&-\displaystyle\int^s_0 (f (t, \sigma^+_i t, W_G (t,
\sigma^+_i t + 0)) -\sigma^+_i W_G (t, \sigma^+_i t + ))   \theta (t, \sigma^+_{i } t ) d t \\
& +\displaystyle O(1) s^3 \, ||\theta||_{{\cal C}^0}
\end{split}
\label{2.16}
\end{eqnarray}
\noindent
for $i = 1, 2,\ldots, p -1$. In the same fashion one can show that
\begin{eqnarray}
\begin{split}
\Delta^1_0 (s, r ; \theta) \;=\;  &\int^{\sigma^-_{1} s}_{- r}W_G (s, x) \theta (s, x) d x
-\int^0_{- r} W_G (0, x) \theta (0, x) d x \\ &+\int^s_0 (f (t,
\sigma^-_{1} t, W_G (t, \sigma^-_1 t - ))
-\sigma^-_1  W_G (t, \sigma^-_1 t  - )) \theta (t, \sigma^-_1 t ) d t\\
&-\int^s_0 f (t, - r, W_G (t, - r)) \theta (t, - r) d t \\
&+ O(1) s^2 (s + r) ||\theta||_{{\cal C}^0}+ O(1) s r^2
||\theta||_{{\cal C}^0},
\end{split}
\label{2.17}
\end{eqnarray}
and
\noindent
\begin{eqnarray}
\begin{split}
\Delta^1_p (s, r ; \theta) \;=\; &\int^r_{\sigma^+_{p} s}W_G (s, x) \theta (s, x) d x -\int^r_0 W_G (0, x) \theta (0, x) dx\\
&+ \int^s_0 f (t, r, W_G (t, r)) \theta (t, r) d t\\
&-\int^s_0 (f (t, \sigma^+_p t, W_G (t,  \sigma ^+_pt + ))
-\sigma^+_p   W_G (t, \sigma^+_p  t + )) \theta(t, \sigma^+_p t) d t\\
&+O(1) s^2(s+r)||\theta||_{{\cal C}^0}+ O(1) s
r^2||\theta||_{{\cal C}^0}.
\end{split}
\label{2.18}
\end{eqnarray}\\

Next, suppose that $W_{C}(t,x)$ consists of an $i$-rarefaction
wave in the region $\big\{(t,x)| \frac{x}{t} \in [\sigma_i^-,
\sigma_i^+]
\big\}$ for some $i \in 1, \ldots, p$. It follows that $W_G(t,x)$ in this region is of the form $$ W_G (t, x) =
\widetilde{W}_C (\frac {x}{t}) + t
\, q (0, 0 ; \widetilde{W}_C (\frac {x}{t})) $$ where
$\widetilde{W}_C (\frac {x}{t})$ is the $i$-rarefaction wave of the
classical Riemann problem $R_C(u_L, u_R ; t_0, x_0)$. By setting
$\xi
=
\frac {x}{t}$, $W_G (t, x) =
\widetilde{W}_G (t, \xi)$, and the technique of change of variables $(t, x) \to
(t, \xi)$, we obtain
\begin{eqnarray}
\begin{split}
&\partial_t W_G + \partial_x f (t, x, W_G)  - g (t, x, W_G)\\
&= \partial_t \widetilde{W}_G - \frac {\xi}{t} \partial_\xi \widetilde{W}_G + \frac {1}{t} \partial_\xi f (t, t \xi, \widetilde{W}_G) - g (t, t \xi, \widetilde{W}_G)\\
&= \frac {1}{t} \Big( \frac {\partial f}{\partial u} (t, t \xi,
\widetilde{W}_G) - \xi I )(I  + t \frac{\partial q}{\partial u}
(0, 0, \widetilde{W}_C)\Big) \cdot
\frac {d \widetilde{W}_C}{d \xi} + q (0, 0, \widetilde{W}_C) -q (t, t \xi, \widetilde{W}_C)
\end{split}
\label{2.19}
\end{eqnarray}
where $I$ is the $p \times p$ identity matrix.
\noindent
Since $\widetilde{W}_C(\xi)$ is a rarefaction wave for the system
\eqref{2.3}, this implies that
\begin{eqnarray}
\frac {1}{t} \Big( -\xi \cdot I+ \frac {\partial f}{\partial u} (0, 0,
\widetilde{W}_C) \Big) \cdot
\frac{d \widetilde{W}_C}{d \xi} = 0.
\label{2.20}
\end{eqnarray}
Thus, by applying \eqref{2.20} to \eqref{2.19} we obtain
\begin{eqnarray}
\begin{split}
\partial_t & W_G + \partial_x f (t, x, W_G) - g (t, x, W_G)\\
=\;&\frac{1}{t} \Big(
\frac{\partial f}{\partial u}(t, x, W_G) - \frac{\partial f}{\partial u}(0, 0, \widetilde{W}_C) \Big)
 \cdot \frac{d \widetilde{W}_C}{d \xi} \\
 &+(\frac{\partial f}{\partial u}(t, x, W_G)-\xi I )
  \frac{\partial q}{\partial u} (0, 0, \widetilde{W}_C)\cdot \frac{d \widetilde{W}_C}{d \xi}
  +q (0, 0, \widetilde{W}_C) - q (t, x, W_G).
\end{split}
\label{2.21}
\end{eqnarray}
\noindent
Next, we multiply \eqref{2.21} by $\theta (t, x)$ and integrate
the equation over the region of $i$-rarefaction wave: $t < s$ and
$\frac{x}{t}
\in [\sigma^-_i,
\sigma^+_i]$. Due to the Lipschitz continuity of
$\frac{\partial f}{\partial u}$ and the fact that $\frac{\partial
f}{\partial u}$, $\frac{\partial q}{\partial u}$, $\frac{d
\widetilde{W}_C}{d \xi}$ and $q$ remain bounded in $[0, s] \times [-r, r]$,
the right hand side of \eqref{2.21} is bounded by $ O(1) s^2 (s +
|u_i - u_{i - 1}|).$
Therefore, by \eqref{2.21} again, we deduce the estimate
\begin{eqnarray}
\begin{split}
\Delta^2_i (s, r ; \theta) =\;& \int^{\sigma^+_i s}_{\sigma^-_i s}W_G (s, x) \theta (s, x) d x \\
&+\int^s_0 (f (t, \sigma^+_i t, W_G (t, \sigma^+_i t))-\sigma^+_i W_G (t, \sigma^+_i t)) \theta (t, \sigma^+_i t) d t\\
&-\int^s_0 (f (t, \sigma^-_i t, W_G (t, \sigma^-_i t)) -\sigma^-_i W_G (t, \sigma^-_i t)) \theta (t, \sigma^-_i t) d t\\
&+O(1) s^2 (s + |u_R - u_L|) \|\theta \|_{{\cal C}^0}.
\end{split}
\label{2.22}
\end{eqnarray}
\noindent
Next, note that an $i$-shock wave satisfies the
Rankine-Hugoniot condition $$ f(0, 0, u_i) -
\sigma_i u_i = f (0,0, u_{i - 1})-\sigma_i u_{i - 1}, $$
\noindent
and this implies that the approximate solution $W_G(t,x)$
satisfies
\begin{eqnarray}
\begin{split}
&\int^s_0  [(f (t, \sigma_i t, W_G (t, \sigma_i t + )) -\sigma_i W_G (t, \sigma_i t + ))] \theta (t, \sigma_i t) d t\\
&-\int^s_0 [(f (t, \sigma_{i-1} t, W_G (t, \sigma_{i-1} t - )) -\sigma_{i-1} W_G (t, \sigma_{i-1} t- ))] \theta (t, \sigma_{i-1} t) d t\\
&=O(1) s^2 |u_R - u_L| \|\theta \|_{{\cal C}^0}
\end{split}
\label{2.23}
\end{eqnarray}
\noindent
where the bound $O(1)$ depends on the Lipschitz constant of $f$
and $L^{\infty}$-norm of $q$. Finally, by the estimates
\eqref{2.12}, \eqref{2.16}-\eqref{2.18} and \eqref{2.22}-\eqref{2.23}, we obtain
\begin{eqnarray*}
\begin{split}
\Delta (s, t ; \theta) =\;& \sum^p_{i = 0} \Delta^1_i (s, t ; \theta) + \sum_{i-rare.\atop waves} \Delta^2_i (s, t ; \theta) \\
=\;&\int^{\sigma^-_i s}_{- r} W_G (s, x) \theta (s, x) d x +\sum^{p - 1}_{i = 1} \int^{\sigma^-_{i + 1} s}_{\sigma^+_i s}W_G (s, x) \theta (s, x) d x\\
&+ \sum_{i-rare.\atop waves} \int^{\sigma^+_i s}_{\sigma^-_i s}
W_G (s, x) \theta (s, x) d x+ \int^r_{\sigma^+_{p} s} W_G (s,
x) \theta (s, x) d x\\
&-\int^0_{- r} W_G (0,x) \theta (0,x)d x -\int^r_0 W_G (0, x) \theta (0, x) d x\\
&+\int^s_0 f (t, r, W_G (t, r)) \theta (t, r) d t - \int^s_0 f (t, - r, W_G (t, - r)) \theta (t, -r) d t\\
&+O(1) (s^2+r^2) (s + r + |u_R - u_L|) \ ||\theta||_{{\cal C}^1},
\end{split}
\end{eqnarray*}
\noindent
which leads to \eqref{2.11} and completes the proof.
\end{proof}
\vskip 0,5cm

%
\vskip 1,0cm
 \setcounter{section}{3}
 \noindent{\bf \large{ 3. Wave interaction estimates}} \\
 \setcounter{equation}{0}

In this section we study the nonlinear interaction of waves issuing from two 
 Riemann solutions and we derive estimates on the wave strengths. 
 
We emphasize that the generalized Riemann solution, nor the 
approximate solution $W_G(t,x)$ of the generalized Riemann problem
$R_G(u_L, u_R ;t_0, x_0)$ is not self-similar. The solution does not consist of regions of
constant value separated by straight lines. We thus should be careful in defining the wave strengths
In fact, we still define here the wave strengths by using the underlying, classical Riemann solution
$W_C(t,x)$. We will see later that this strategy is accurate enough and that 
the discrepancy in total variation between $W_G(t,x)$ and $W_C(t,x)$ on each time
step is uniformly small (Cf.~Section 4) 
when our Glimm scheme is applied to the problem \eqref{1.1}, \eqref{1.2}. 
The same observation applies to the potential of wave interaction to be introduced later. 

In the rest of the section, all waves are considered as waves from some classical
Riemann problem unless specified otherwise. We say that an $i$-wave and a
$j$-wave approach each other (or interact in the future) if either $i > j$, or else $i = j$ and
at least one of two waves is a shock wave. Suppose there are two
solutions from different classical Riemann problems with strengths denoted by $\alpha = (\alpha_i, \ldots,
\alpha_p)$ and $\beta = (\beta_i, \ldots, \beta_p)$, then
the {\sl wave interaction potential}  associated 
these two solutions is defined by
\begin{eqnarray}
D (\alpha, \beta) := \sum_{(i,j)} |\alpha_i
\beta_j|, 
\label{3.1}
\end{eqnarray}
\noindent
where the notation $(i,j)$ under the summation sign indicates an
$i$-wave in one solution approaching a $j$-wave in the other
solution, and the summation is on all approaching waves; also
$\alpha_i$ or $\beta_i$ is negative when $i=j$. In addition,
given a $(u_L, u_R ; t_0, x_0) \in {\cal U} \times {\cal U}
\times \RR_+ \times \RR$, the wave strengths in $R_C (u_L, u_R ; t_0, x_0)$
are denoted by $\eps (u_L, u_R ; t_0,
x_0)$.

We first recall:

\vskip0,5cm
\noindent
{\bf Lemma 3.1.} (Glimm) {\em \quad 1) Given a $(t_0, x_0)$ in
$\RR_+ \times
\RR$ and $u_L$, $u_M$, $u_R$ in ${\cal U}$,
we have
\begin{eqnarray}
|\gamma - (\alpha + \beta)|= O(1) D (\alpha, \beta)
\label{3.2}
\end{eqnarray}
where
\begin{eqnarray}
\alpha = \eps (u_L, u_M ; t_0, x_0), \quad
\beta = \eps (u_M, u_R ; t_0, x_0), \quad
\gamma = \eps (u_L, u_R ; t_0, x_0).
\label{3.3}
\end{eqnarray}
2) Let $v_L, v_R$ be two constant states in ${\cal U}$, then
\begin{eqnarray}
D (\gamma, \delta) = D (\alpha, \delta) + D (\beta, \delta) + O(1)
|\delta| D (\alpha, \beta),
\label{3.4}
\end{eqnarray}
and
\begin{eqnarray*}
D (\delta, \gamma) = D (\delta, \alpha) + D (\delta, \beta) + O(1)
|\delta| D (\alpha, \beta)
\label{3.5}
\end{eqnarray*}
where $\alpha, \beta$ and $\gamma$ are given in \eqref{3.3}, and
$\delta$ is given by $\delta = \eps (v_L, v_R ; t_0, x_0)$.
}

The following lemma describes the dependence of the
wave strengths and potential $D(\cdot, \cdot)$ with respect to
their arguments. We introduce the following ``local norm" of a given function
$\varphi(t,x,u)$
\begin{eqnarray}
N^{x_1,x_2}_{t_1,t_2} (\varphi) =
\sup\!\!\!\!\!\!& \Big\{|\varphi (t, x, u)| \; ; t \in  [t_1, t_2],
\; x \in [x_1, x_2],\; u \in {\cal U} \Big\},
\label{3.7}
\end{eqnarray}
\noindent
where the supremum is taken over any function $u \in {\cal U} $ and
$(t,x) \in [t_1, t_2] \times [x_1, x_2]$.
\vskip0,5cm
\noindent

{\bf Lemma 3.2.} {\em
\quad 1) The wave strength $\eps = (\eps_i)_{1
\leq i \leq p} : {\cal U} \times {\cal U} \times
\RR_+ \times \RR \to \RR^p$ is a ${\cal C}^2$ vector function of
its arguments. Furthermore, for any $(u_L, u_R), (u^\prime_L,
u^\prime_R)$ in ${\cal U} \times {\cal U}$ and any $(t_0, x_0),
(t^\prime_0, x^\prime_0)$ in $\RR_+ \times
\RR$, we have
\begin{eqnarray}
\begin{split}
|\alpha^\prime-\alpha | \;=\;\;& O(1) |\alpha| (|u^\prime_L -
u_L| + |u^\prime_R - u_R| + C^0_1 |t^\prime_0 - t_0| + C^0_2
|x^\prime_0 - x_0|)\\
&+ O(1) |(u^\prime_R - u^\prime_L) - (u_R - u_L)|
\end{split}
\label{3.8}
\end{eqnarray}
where
\begin{eqnarray}
\alpha = \eps (u_L, u_R ; t_0, x_0), \quad
\alpha^\prime = \eps (u^\prime_L,
u^\prime_R ; t^\prime_0, x^\prime_0),
\label{3.9}
\end{eqnarray}
and the constants $C^0_1$ and $C^0_2$ are given by
\begin{eqnarray}
C^0_1 := N^{x_0, x_0^\prime}_{t_0, t_0^\prime} (\frac
{\partial^2 A}{\partial t \partial u}), \qquad C^0_2 :=
N^{x_0, x_0^\prime}_{t_0, t_0^\prime} (\frac {\partial^2
A}{\partial x
\partial u}).
\label{3.10}
\end{eqnarray}
\vskip 0,5cm
2) For given $(u_L, u_R), (v_L, v_R), (u^\prime_L, u^\prime_R)$,
$(v^\prime_L, v^\prime_R)$ in ${\cal U} \times {\cal U}$ and
$(t_1, x_1), (t_2, x_2)$, $(t^\prime_1, x^\prime_1)$,
$(t^\prime_2, x^\prime_2)$ in $\RR_+
\times \RR$, we have
\begin{eqnarray}
\begin{split}
D (\alpha^\prime, \beta^\prime) =\;& \displaystyle D(\alpha,
\beta)
+O(1)|\alpha| |(v^\prime_R - v^\prime_L) - (v_R - v_L)|\\
&+O(1) |\beta| |(u^\prime_R - u_L^\prime) - (u_R - u_L)|\\
&+O(1) |\alpha| |\beta|  \Big( |u^\prime_L -u_L| +|u^\prime_R - u_R| + |v^\prime_L - v_L| + |v^\prime_R - v_R|\Big)
\\
&+O(1) |\alpha| |\beta| \sum_{m=1, 2}
\{C^m_1|t^\prime_m - t_m| + C^m_2 |x^\prime_m - x_m|\}\\
&+ O(1)|(u^\prime_R - u^\prime_L) - (u_R - u_L)|\cdot|(v^\prime_R
- v^\prime_L) - (v_R - v_L)|
\end{split}
\label{3.11}
\end{eqnarray}
where
\begin{eqnarray}
\begin{split}
&\alpha = \eps (u_L, u_R ; t_1, x_1), \; \beta = \eps(v_L, v_R ; t_2, x_2),\\
&\alpha^\prime = \displaystyle \eps (u^\prime_L,
u^\prime_R ; t^\prime_1, x^\prime_1), \; \beta^\prime =
\eps (v^\prime_L, v^\prime_R ; t^\prime_2, x^\prime_2),
\end{split}
\label{3.12}
\end{eqnarray}
\noindent
and the constants $C^m_1$, $C^m_2$ are defined by
\begin{eqnarray}
C^m_1 := N^{x_m, x^\prime_m}_{t_m, t^\prime_m}(\frac
{\partial^2 A}{\partial t \partial u}), \quad C^m_2 :=
N^{x_m, x^\prime_m}_{t_m, t^\prime_m}(\frac{\partial^2
A}{\partial x
\partial u}), \ m = 1, 2.
\label{3.13}
\end{eqnarray}
}
\noindent
\begin{proof}
\quad The regularity of functions $\eps_i$, $i = 1, 2, \ldots,
p$, is a consequence of smoothness of the flux function $f$ and
the result of \cite{La2}. Moreover, the functions $\frac {\partial^2
\eps_i}{\partial t \partial u_R}$ and $\frac {\partial^2
\eps_i}{\partial x \partial u_R}$ are bounded if
$\frac {\partial ^2 A}{\partial t \partial u}$, $\frac {\partial^2
A}{\partial x \partial u}$ are bounded.

To show \eqref{3.8}, we
note that $\eps_i (u_L, u_R ; t_0, x_0) = 0$  when $u_R = u_L$.
Then, by the regularity of  $\eps_i$, $i = 1, 2,
\ldots, p$, we can express $\eps_i (u_L, u_R ; t_0, x_0)$,
$\eps_i (u^\prime_L, u^\prime_R ; t^\prime_0, x^\prime_0)$
as
\begin{eqnarray*}
\eps_i (u_L, u_R ; t_0, x_0) = \int^1_0 \frac{\partial \eps_i}{\partial u_R}
(u_L, (1 - \tau) u_L + \tau u_R ; t_0, x_0) d \tau \cdot(u_R-u_L),
\end{eqnarray*}
\begin{eqnarray*}
\eps_i (u^\prime_L, u^\prime_R ; t^\prime_0,
x^\prime_0) =
\int^1_0 \frac{\partial \eps_i}{\partial u_R}
(u^\prime_L, (1 -\tau) u^\prime_L + \tau u^\prime_R ; t^\prime_0,
x^\prime_0) d \tau \cdot(u^\prime_R - u^\prime_L).
\end{eqnarray*}
 Applying the definition of $\{C^0_j: j=1,2.\}$ in \eqref{3.10}
 and the norm in \eqref{3.7}, we obtain
\begin{eqnarray*}
\begin{split}
&\quad \;\;\eps_i (u^\prime_L, u^\prime_R ; t^\prime_0, x^\prime_0)
-\eps_i (u_L, u_R ; t_0, x_0) \\
&=\int^1_0 (\frac {\partial \eps_i}{\partial u_R}
(u^\prime_L, (1 - \tau) u^\prime_L + \tau u^\prime_R ;
t^\prime_0, x^\prime_0) -\frac {\partial \eps_i}{\partial
u_R}
(u_L, (1 - \tau) u_L + \tau u_R ; t_0, x_0)) d \tau \cdot (u_R - u_L) \\
&\quad+\int^1_0 \frac {\partial \eps_i}{\partial
u_R}(u^\prime_L, (1 - \tau) u^\prime_L + \tau u^\prime_R ;
t^\prime_0, x^\prime_0)
d \tau \cdot ((u^\prime_R - u^\prime_L) - (u_R - u_L)) \\
&=O(1) \{|u^\prime_L - u_L| + |u^\prime_R - u_R| +C^0_1 |t^\prime_0 - t_0| + C^0_2 |x^\prime_0 - x_0|\} |u_R - u_L|\\
&\quad+O(1) |(u^\prime_R - u^\prime_L) - (u_R - u_L)|,
\end{split}
\end{eqnarray*}
\noindent
Therefore, by the observation of \eqref{3.9} and the fact that
\begin{eqnarray*}
|u_R- u_L| = O(1) |\eps (u_L, u_R ; t_0, x_0)| = O(1)
|\alpha|,
\end{eqnarray*}
we obtain \eqref{3.8}.\\

Next we derive \eqref{3.11}. By applying \eqref{3.8} directly, we have
\begin{eqnarray*}
\begin{split}
\alpha^\prime_i \;=\;&\alpha_i + O(1) |\alpha| \{|u^\prime_L - u_L| +
|u^\prime_R - u_R| + C^1_1 |t^\prime_1 - t_1|+ C^1_2 |x^\prime_1 - x_1|\}\\
&+O(1) |(u^\prime_R - u^\prime_L) - (u_R - u_L)|,
\end{split}
\end{eqnarray*}
\noindent

\begin{eqnarray*}
\begin{split}
\beta^\prime_j \;=\;& \beta_j + O(1) |\beta|
\{|v^\prime_L - v_L| + |v^\prime_R - v_R| + C^2_1 |t^\prime_2 - t_2|
+ C^2_2 |x^\prime_2 - x_2|\}\\
&+O(1)|(v^\prime_R -
v^\prime_L) - (v_R - v_L)|
\end{split}
\end{eqnarray*}
\noindent
for $i,\; j = 1, 2,\ldots, p$ where the constants $\{ C^m_j: j,
m=1,2.\} $ are given in \eqref{3.13} and \eqref{3.7}. We define $ A
:= \{|u^\prime_L - u_L| + |u^\prime_R - u_R|+ C^1_1
|t^\prime_1 - t_1| + C^1_2 |x^\prime_1 - x_1| \}$ and $B :=
\{|v^\prime_L - v_L| + |v^\prime_R - v_R| +C^2_1 |t^\prime_2 -
t_2| + C^2_2 |x^\prime_2 - x_2|\}$. Then by multiplying two
previous equations together and using the fact that $A, B$ are of
order $O(1)$ for $(u_L, u_R), (u^\prime_L, u^\prime_R) \in
\cal U \times \cal U$, we obtain
\begin{eqnarray*}
\begin{split}
\alpha^\prime_i \beta^\prime_j
\;=\; & \alpha_i \beta_j + O(1)|\alpha||\beta|(A+B)+O(1)|\alpha||(v^\prime_R - v^\prime_L) -
(v_R - v_L)|\\
&+O(1)|\beta||(u^\prime_R - u^\prime_L) - (u_R -
u_L)|\\
&+O(1)|(u^\prime_R - u^\prime_L) - (u_R - u_L)|\cdot|(v^\prime_R
- v^\prime_L) - (v_R - v_L)|,
\end{split}
\end{eqnarray*}
\noindent
$i, j= 1, 2, \ldots, p$. Summing up previous equations for $i,\;
j = 1, 2,\ldots, p$, we obtain \eqref{3.11}. The proof is
completed.
\end{proof}

Using Lemmas 3.1 and 3.2, we obtain wave interaction estimates --which can be interpreted as
a generalized version of \cite{Gl}.
\vskip 0,5cm
\noindent
{\bf Proposition 3.3. }{\em \quad  1) Suppose that $s$, $r$ are two
positive numbers and $(t_0, x_0)$ is in $\RR_+ \times
\RR$. Also assume that $u_L$, $u_M$, $u_R$, $u_L +
\mu_L$, $u_R + \mu_R $ are constant states in
${\cal U}$ and $\alpha,
\beta$ and $\gamma$ are the wave strengths of solutions of three
classical Riemann problems $R_C (u_L, u_M ; t_0, x_0 - r)$, $R_C
(u_M, u_R ; t_0, x_0 + r)$ and $R_C (u_L + \mu_L, u_R +
\mu_R ; t_0 + s, x_0)$, i.e.,
\begin{eqnarray}
\begin{split}
&\alpha = \eps (u_L, u_M ; t_0, x_0 - r), \quad
\beta =  \eps (u_M, u_R ; t_0, x_0 + r),\\
&\gamma = \eps (u_L + \mu_L, u_R + \mu_R ; t_0 + s,
x_0).
\end{split}
\label{3.14}
\end{eqnarray}
\noindent
Then we have
\begin{eqnarray}
\begin{split}
|\gamma|\; =\; &|\alpha| + |\beta| + O(1) D (\alpha, \beta)\\
\;\;\; &+ O(1) (|\alpha| + |\beta|) (|\mu_L| + |\mu_R| +C_1 s + C_2 r)\\
&+ O(1)|\mu_R - \mu_L|
\end{split}
\label{3.15}
\end{eqnarray}
\noindent
where constants $C_1$ and $C_2$ are defined by
\begin{eqnarray}
C_1 := N^{x_0 x_0}_{t_0, t_0 + s} \bigl( \frac {\partial^2
A}{\partial t
\partial u}\bigr), \quad C_2 := N^{x_0 - r, x_0 + r}_{t_0, t_0} \bigl(\frac
{\partial^2 A}{\partial x \partial u}\bigr).
\label{3.16}
\end{eqnarray}
\noindent
 2) Let $\alpha$, $\beta$, $\gamma$ be the wave strengths as described in \eqref{3.14}.
 Also, for a given $(v_L, v_R)$ in ${\cal U} \times
{\cal U}$ and $(t_1, x_1)$ in $\RR_+ \times
\RR$, we define $ \delta = \eps (v_L, v_R ; t_1, x_1).$ Then
\begin{eqnarray}
\begin{split}
D (\gamma, \delta)= \; &  \displaystyle D (\alpha, \delta)  + D (\beta, \delta) + O(1)|\delta| D (\alpha, \beta) +O(1) |\delta| |\mu_R - \mu_L|\\
&+\displaystyle O(1) |\delta|(|\alpha| + |\beta|) (|\mu_L| + |\mu_R| + C_1 s + C_2 r),\\
\end{split}
\label{3.17}
\end{eqnarray}
and
\begin{eqnarray}
\begin{split}
D (\delta, \gamma)= \; &  \displaystyle D (\delta, \alpha)  + D (\delta, \beta) + O(1)|\delta| D (\alpha, \beta)
+O(1) |\delta| |\mu_R - \mu_L|\\
&+\displaystyle O(1) |\delta|(|\alpha| + |\beta|) (|\mu_L| + |\mu_R| + C_1 s + C_2 r)\\
\end{split}
\label{3.18}
\end{eqnarray}
\noindent
where constants $C_1$ and $C_2$ are given in \eqref{3.16}. }
\vskip0,5cm
\begin{proof}
\quad By the definition of $\gamma$ in \eqref{3.14} and Lemma 3.2
 with $u^\prime_L = u_L + \mu_L$, \\ $u^\prime_R = u_R +
\mu_R, \ t^\prime = t_0 + s, \ x^\prime_0 = x_0$, we obtain
\begin{eqnarray}
\begin{split}
\gamma =\; & \;  \eps (u_L, u_R ; t_0, x_0) + O(1) |\eps
(u_L, u_R ; t_0, x_0)| \, \Big( |\mu_L| + |\mu_R| + C_1 \, s\Big)
\\
 &+O(1) \, |\mu_R - \mu_L|
\end{split}
\label{3.19}
\end{eqnarray}
\noindent
where constant $C_1$ is given in \eqref{3.16}. Similarly, by Lemma 3.2
we have
\begin{eqnarray}
\eps(u_L, u_M ; t_0, x_0)= \alpha  + O(1) C_2 |\alpha| r,
\label{3.20}
\end{eqnarray}
\begin{eqnarray}
\eps (u_M, u_R ; t_0, x_0)= \beta + O(1) C_2 |\beta| r.
\label{3.21}
\end{eqnarray}
\noindent
On the other hand, Glimm's interaction estimates \eqref{3.2}, \eqref{3.3} lead to
\begin{eqnarray}
\begin{split}
\eps (u_L, u_R ; t_0, x_0) = \; & \; \displaystyle  \eps (u_L, u_M ; t_0, x_0) + \eps (u_M, u_R ; t_0, x_0)\\
&+\displaystyle O(1) D (\eps (u_L, u_M ; t_0, x_0),\eps (u_M, u_R ; t_0, x_0)).
\end{split}
\label{3.22}
\end{eqnarray}
\noindent
Also, by \eqref{3.11}-\eqref{3.13} with $\alpha^\prime =
\eps (u_L, u_M ; t_0, x_0)$ and $ \\
 \beta^\prime = \eps (u_M, u_R ; t_0, x_0)$, we obtain
\begin{eqnarray}
D (\eps (u_L, u_M ; t_0, x_0),
\eps (u_M, u_R ; t_0, x_0)) = D (\alpha, \beta) +
O(1) |\alpha| |\beta| C_2 r.
\label{3.23}
\end{eqnarray}
\noindent
Then, from \eqref{3.19}-\eqref{3.23} it follows that
\begin{eqnarray} 
\begin{split}
|\gamma| \; &= \; \displaystyle |\alpha| + |\beta| + O(1)D (\alpha, \beta) + O(1) (|\alpha|+|\beta|+|\alpha| |\beta|)
 C_2 r \\   
&\quad \  +\displaystyle O(1) |\eps (u_L, u_R ; t_0, x_0)| (|\mu_L| + |\mu_R| + C_1 s)
       + \displaystyle O(1) |\mu_R-\mu_L|\\
\; &= \; \displaystyle |\alpha| + |\beta| + O(1)D (\alpha,
\beta) + O(1) (|\alpha|+|\beta|)
 C_2 r \\   
&\quad \  +\displaystyle O(1) |\eps (u_L, u_R ; t_0, x_0)| (|\mu_L| + |\mu_R| + C_1 s)
  + \displaystyle O(1) |\mu_R-\mu_L|.
\end{split}
\label{3.24}
\end{eqnarray}
\noindent
Also, we see that estimates \eqref{3.22} and \eqref{3.23} yield
\begin{eqnarray*}
\begin{split}
|\eps (u_L, u_R ; t_0, x_0)| \;=& \; |\eps (u_L, u_M ;t_0,x_0)| + |\eps (u_M, u_R ; t_0, x_0)|\\
&+ O(1) D(\eps (u_L, u_M ; t_0, x_0),
\eps(u_M,u_R; t_0, x_0)) \\ =& |\alpha| + |\beta| + O(1) D
(\alpha,\beta) + O(1) \, (|\alpha| + |\beta|) \, C_2 r \\
&+ O(1) |\alpha| |\beta| C_2 r\\ =& (|\alpha| + |\beta|) (1 +
O(1) C_2 r) + O(1)D(\alpha,\beta)+ O(1) |\alpha| |\beta|C_2 r,
\end{split}
\end{eqnarray*}
\noindent
which in particular implies that
\begin{eqnarray}
|\eps (u_L, u_R ; t_0, x_0)| = O(1) (|\alpha| + |\beta|).
\label{3.25}
\end{eqnarray}
Therefore, combining \eqref{3.24} with \eqref{3.25}, we obtain
\eqref{3.15}.\\

Next we derive \eqref{3.17}. The proof of \eqref{3.18} is similar, and is omitted.
By the estimate \eqref{3.4} we see that
\begin{eqnarray}
\begin{split}
D(\eps (u_L, u_R ; t_0, x_0), \delta) \;=& \; \displaystyle D (\eps (u_L, u_M ; t_0, x_0), \delta)
+D (\eps (u_M, u_R ; t_0, x_0), \delta) \\
 &+ \displaystyle O(1) |\delta| D (\eps (u_L, u_M ; t_0, x_0), \eps (u_M, u_R ; t_0, x_0)).
\end{split}
\label{3.26}
\end{eqnarray}
On the other hand, estimate \eqref{3.11} yields
\begin{eqnarray}
\begin{split}
D (\gamma, \delta)\;=& \; \displaystyle D(\eps (u_L, u_R ; t_0, x_0), \delta) + O(1) |\delta| |\mu_R - \mu_L|\\   
&+ \displaystyle O(1) |\eps (u_L, u_R ; t_0, x_0)|
|\delta| (|\mu_L| + |\mu_R| + C_1 s),
\end{split}
\label{3.27}
\end{eqnarray}
\begin{eqnarray}
D (\eps (u_L, u_M ; t_0, x_0), \delta) = D (\alpha, \delta)
+ O(1) |\alpha| |\delta| C_2 r,
\label{3.28}
\end{eqnarray}
\begin{eqnarray}
D (\eps (u_M, u_R ; t_0, x_0), \delta) = D (\beta, \delta)
+ O(1) |\beta| |\delta| C_2 r,
\label{3.29}
\end{eqnarray}
and
\begin{eqnarray}
D (\eps (u_L, u_M ; t_0, x_0), \eps (u_M, u_R ;
t_0, x_0)) = D (\alpha, \beta) + O(1) |\alpha| |\beta| C_2 r.
\label{3.30}
\end{eqnarray}
Thus, by applying \eqref{3.25}, \eqref{3.27}-\eqref{3.30} to
\eqref{3.26}, we obtain the estimate \eqref{3.17}. The proof is
completed.
\end{proof}
\quad We just showed in Proposition 3.3 that Glimm's interaction estimates
(Lemma 3.1) remain valid for the quasilinear hyperbolic system
\eqref{1.1} up to certain error terms. The following 
immediate consequence of Proposition 3.3 will be the key to the forthcoming stability result.
\vskip0,5cm
\noindent
{\bf Corollary 3.4} {\em
\quad Following the notations and assumptions in Proposition 3.3 and letting
\begin{eqnarray*}
\mu_L := - s q (t_0 + s, x_0, u_L), \ \mu_R := - s q
(t_0 + s, x_0, u_R)
\label{3.31}
\end{eqnarray*}
in (3.14), we have
\begin{eqnarray}
\begin{split}
|\gamma|\; = \;\ &|\alpha| + |\beta| + O(1) D (\alpha, \beta) \\
&+O(1) (|\alpha| + |\beta|) \Big( (C_1 + C_3 + C_4) \, s + C_2 \, r\Big),
\end{split}
\label{3.32}
\end{eqnarray}
\begin{eqnarray}
\begin{split}
D (\delta, \gamma)\;= \; \ &D (\delta, \alpha)  + D (\delta, \beta)+ O(1)|\delta| D (\alpha, \beta) \\
&+ O(1) |\delta|(|\alpha|+|\beta|) \Big( (C_1+C_3+C_4) \, s+C_2 \, r \Big),
\end{split}
\label{3.33}
\end{eqnarray}
\begin{eqnarray}
\begin{split}
D (\gamma, \delta)\;=\; \ & D (\alpha, \delta)  + D (\beta, \delta) + O(1)|\delta| D (\alpha, \beta) \\
&+ O(1) |\delta|(|\alpha|+|\beta|) \{ (C_1+C_3+C_4)s+C_2r \}
\end{split}
\label{3.34}
\end{eqnarray}
\noindent
where constants $C_1$, $C_2$ are given in \eqref{3.16} and $C_3$,
$C_4$ are given by
\begin{eqnarray}
C_3 := N^{x_0 x_0}_{t_0, t_0 + s} (q(t,x,u)), \quad C_4
:= N^{x_0 x_0}_{t_0, t_0 + s}(\frac {\partial q}{\partial
u}(t,x,u)).
\label{3.35}
\end{eqnarray}
\begin{proof}
By the observation of \eqref{3.25} we obtain
\begin{eqnarray}
\begin{split}
|\mu_R-\mu_L| \;=\ & s|q(t_0+s,x_0,u_R)-q(t_0+s,x_0,u_L)|\\
=\ & s\frac {\partial q}{\partial u}(t_0+s,x_0,\bar{u})\cdot
|u_R-u_L|\\
=\ & O(1)C_4 s|\eps(u_L,u_R;t_0,x_0)|\\
=\ & O(1)(|\alpha|+|\beta|)C_4 s
\end{split}
\label{3.36}
\end{eqnarray}
where $\bar{u} \in \cal U$ and $C_4$ is given in \eqref{3.35}.
Therefore, by combining \eqref{3.36} with the result of Proposition 3.3, we obtain \eqref{3.32}-\eqref{3.34} . The proof is completed.
\end{proof}
}
\vskip0,5cm
%
%
 \setcounter{section}{4}
 \noindent{\bf \large{4. Stability of the generalized Glimm method}} \\
 \setcounter{equation}{0}

We are in position to introduce our version of Glimm scheme for the
approximation of the quasilinear system \eqref{1.1}. Then we rely on 
the wave interaction estimates in Section 3 and prove a stability result.
\vskip0,5cm

The approximate solution to the Cauchy
problem \eqref{1.1}, \eqref{1.2} is defined as follows.
Given two positive constants $s$ and $r$ satisfying the C-F-L
condition \eqref{2.10}, we introduce the constant
\begin{eqnarray}
\lambda_* := \frac{r}{s}.
\label{4.1}
\end{eqnarray}
\noindent
Let also $a = \big\{a_k: a_k \in (- 1, 1),k \in \mathbb{N} \big\}$
 be an equidistributed sequence.
We divide the $(t,x)$ plane into
\begin{eqnarray}
t_k=ks, \; x_h=hr, \quad k=0,1,2, \ldots, \; h \in \mathbb{Z}.
\label{4.2}
\end{eqnarray}
Next, we construct an approximate solution $u_r (t, x)$ of the
problem \eqref{1.1}, \eqref{1.2} in the following way. First,
the initial data $u_0(x)$ is approximated by a piecewise constant
function
\begin{eqnarray}
u_r (0, x) = u_0 (h r), \quad x \in [(h-1)r,(h + 1)r ),
\; h \; \rm is
\; \rm odd.
\label{4.3}
\end{eqnarray}
Then, within domain $0 \leq t < s$, we construct an approximate
solution $W_G(t,x)$ for each generalized Riemann problem with
initial data $u_r(0,x)$ to obtain $u_r(t,x)$ in region $\big\{(t,x);
0 \leq t < s
\big\}$. If $u_r (t, x)$ has been constructed for $t < k s$, $k
\in
\mathbb{N}$, we set
\begin{eqnarray}
\begin{split}
u_r (k s, x) \;:= \;& u_r (k s -, (h + a_k) r)
\end{split}
\label{4.4}
\end{eqnarray}
for $x \in [(h - 1) r, (h + 1) r)$, $k + h $ is odd. Again, we
solve the generalized Riemann problems with initial data $u_r(ks,
x)$ given in \eqref{4.4} to construct $u_r(t,x)$ within region $\big\{(t,x);
ks
\leq t < (k+1)s \}$. Following the process \eqref{4.3}, \eqref{4.4} consecutively, we then construct
our approximate solution $u_r(t,x)$ of \eqref{1.1}, \eqref{1.2}. In other
words, the approximate solution
to the problem \eqref{1.1}, \eqref{1.2} generated by the generalized Glimm scheme is given by
\begin{eqnarray}
\begin{split}
u_r (t, x)=\;& W_G (t, x ; u_r (k s, (h - 1) r),u_r (k s, (h + 1)
r) ; k s, h r)
\end{split}
\label{4.5}
\end{eqnarray}
for $(t,x)
\in [k s, (k + 1) s) \times [(h - 1) r,(h + 1) r)$,
$k + h $ is even.\\

Next we study the stability of $u_r(t,x)$ in $L^{\infty}$ and $BV$
norms. This requires the description of mesh points, mesh curves
and immediate successors beforehand. Recall that the values of
$u_r(t,x)$ on $t=ks$ are determined by the values of $u_r(t,x)$ at
points $\{(k s-, (h + a_k) r); \; h \in \mathbb{Z}, \; k+h \;
\rm is \; \rm odd \}$, we call these points $\{(k s, (h + a_k) r): k=0,1,2,
\cdots,
\; h \in \mathbb{Z}, \; k+h \; \rm is \; \rm odd
\}$ the $mesh \; points$ of approximate
solution $u_r(t,x)$. We obtain a set of diamond regions by
connecting all mesh points with segments. An unbounded piecewise
linear curve $I$ is called a $mesh
\; curve$ if $I$ lies on the boundaries of those diamond regions.
Suppose $I$ is a mesh curve, then $I$ divides the $(t,x)$ plane
into $I^{+}$ and $I^{-}$ regions, such that $I^{-}$ contains
$t=0$. We say two mesh curves $I_1 > I_2$ ( $I_1$ is a $successor$
of $I_2$) if every point of $I_1$ is either on $I_2$ or contained
in $I^{+}_2$. And, $I_1$ is an $immediate \; successor$ of $I_2$
if $I_1 > I_2$ and every mesh point of $I_1$ except one is on
$I_2$. Note that the difference between $I_1$ and $I_2$ is
determined by a diamond region if one
is an immediate successor of the other. \\

Next, to simplify the notations, we set $u_{k, h} :=  u_r (k s, h r)$ when $k + h$ is odd.
By the observation of \eqref{2.7} and \eqref{4.5}, we have 
\begin{eqnarray*}
u_{k, h} = \tilde{u}_{k, h} + s \, q ((k - 1) s, h r,
\tilde{u}_{k, h}), \quad k + h \ \rm is \ \rm odd,
\label{4.7}
\end{eqnarray*}
where $\tilde{u}_{k, h}$ is the value of $R_C (u_{k - 1, h - 1},
u_{k - 1, h+1} ; (k - 1) s, h r)$ at $(ks-,(h+a_k)r)$, i.e.,
\begin{eqnarray*}
\tilde{u}_{k, h} = W_C (a_k \frac {r}{s} ; u_{k - 1, h - 1},
u_{k - 1, h+1} ; (k - 1) s, h r)
\label{4.8}
\end{eqnarray*}
with the function $W_C$ given in Section 2. Next, given a pair
$(k_0, h_0)$, $k_0+h_0 \; \rm is \; \rm even$, we note that the
$(t, x) -$plan consists of the diamond regions $\Gamma _{k_0,
h_0}$ with center $(k_0 s, h_0 r)$ and vertices (mesh points)
\begin{eqnarray}
\begin{split}
&S := ((k_0 - 1) s,(h_0 + a_{k_0 - 1}) r),\;W
:= (k_0 s,(h_0 - 1 + a_{k_0}) r),  \\
 &E := (k_0 s,(h_0 + 1 + a_{k_0}) r), \;N
:= ((k_0 + 1) s,(h_0 + a_{k_0 + 1}) r)
\end{split}
\label{4.9}
\end{eqnarray}
(see Figure 4.1). We set
\begin{eqnarray}
u_S := u_{k_0 - 1, h_0}, \ u_W := u_{k_0, h_0 - 1},
\ u_E := u_{k_0, h_0 + 1}, \ u_N := u_{k_0 + 1, h_0},
\label{4.10}
\end{eqnarray}
and
\begin{eqnarray}
\tilde{u}_S := \tilde{u}_{k_0-1, h_0 }, \
\tilde{u}_W := \tilde{u}_{k_0, h_0 - 1}, \
\tilde{u}_E := \tilde{u}_{k_0, h_0 + 1}, \
\tilde{u}_N := \tilde{u}_{k_0 + 1, h_0}.
\label{4.11}
\end{eqnarray}
Note that $u_W$ and $u_E$ are the states in $R_G((k_0 - 1) s,
(h_0 - 1) r)$ and $R_G((k_0 - 1) s, (h_0 + 1) r)$ respectively,
i.e.,
\begin{eqnarray*}
u_W = \tilde{u}_W + s \, q ((k_0 - 1) s, (h_0 - 1) r, \tilde
u_W),
\label{4.12}
\end{eqnarray*}
\begin{eqnarray*}
u_E = \tilde{u}_E + s \, q ((k_0 -1) s, (h_0 + 1) r,
\tilde{u}_E).
\label{4.13}
\end{eqnarray*}
\begin{center}
\includegraphics[height = 5cm, width=12 cm]{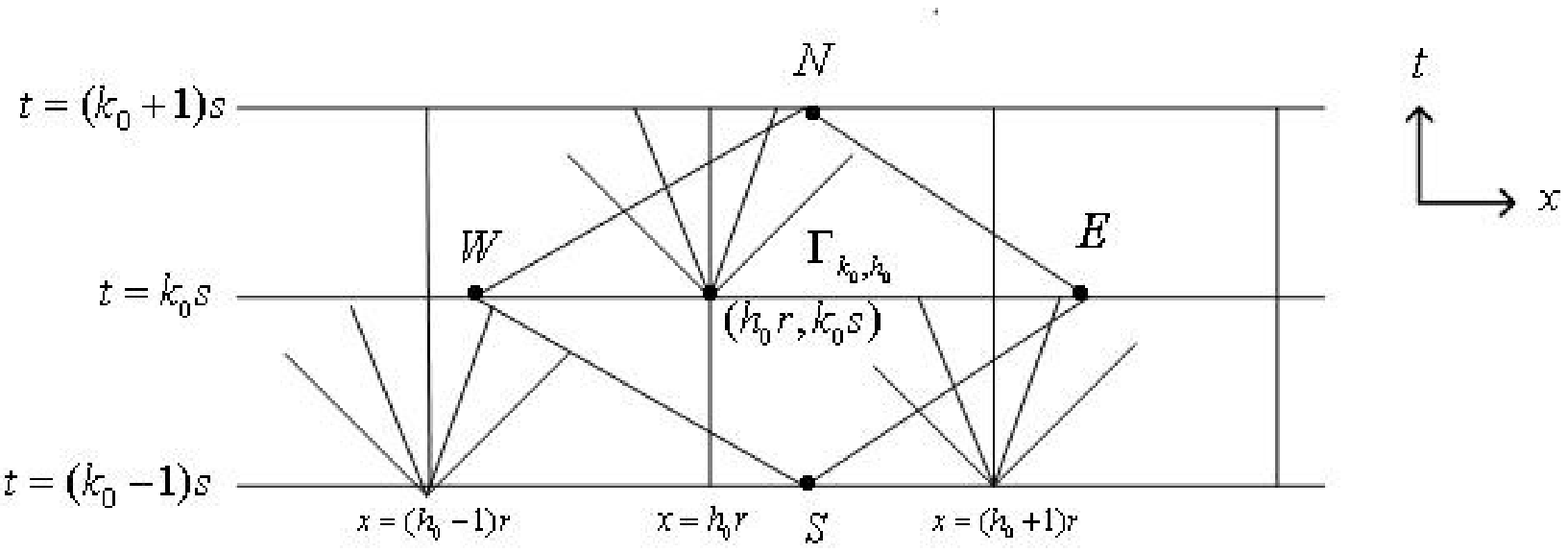} \\
  Figure 4.1 : Diamond region $\Gamma_{k_{0},h_{0}}$
\end{center}

Now we define the $strengths$ of waves in $u_r(t,x)$. However, the
set up for the waves strengths of $u_r(t,x)$ becomes crucial due
to the lack of self-similarity of approximate solution
$W_G(t,x)$, the strengths of waves in $W_G(t,x)$ can not be
defined in the traditional way as described in \cite{La2}. To
overcome the difficulty, we first solve the associated classical
Riemann problems with the initial data $\{u_r (k s -, (h + a_k)
r); \; x
\in [(h - 1) r, (h + 1) r),
\; k+h \ \rm is  \ \rm odd \}$ (see \eqref{4.4}) within each
time step. So we construct a new function $\tilde{u}_r(t,x)$
defined on $\RR_+ \times
\RR$. Then we define the strengths of
approximate waves in $u_r(t,x)$ based on classical waves in
$\tilde{u}_r(t,x)$. More precisely, given a wave $(u_{i-1}(t),
u_i(t))$ in $u_r(t,x)$, there exist two corresponding constant
states $u_{i-1}$, $u_i$ and a classical Riemann wave $(u_{i-1},
u_i)$ with strength $\eps (u_{i-1}, u_i)$ in $\tilde{u}_r(t,x)$,
then the strength of $(u_{i-1}(t), u_i(t))$ is defined as
$\eps (u_{i-1}, u_i)$. \\

Next, we show that, under the condition that the $L^1(\mathbb{R_+}
\times \RR)$-norms of $q$ and $\frac{\partial{q}}{\partial{u}}$
are small, the sum of strengths for waves in $u_r(t,x)$ crossing
mesh curve $J$ can be regarded as an equivalent norm for the total
variation of $u_r(t,x)$ on $J$. By the fact that the term
$|\eps (u_{i-1}, u_i)|$ is equivalent to the total
variation of $(u_{i-1}, u_i)$ for any classical Riemann wave
$(u_{i-1}, u_i)$, it is equivalent to show that the total
variation of $u_r(t,x)$ on $J$ is equivalent to the total
variation of $\tilde{u}_r(t,x)$ on $J$. To show this, let $J_k$
be a mesh curve lying within $k$-th time level $\{(t,x); \; ks
\leq t < (k+1)s \}$, and let $TV(u_r(t,x), J_k)$,
$TV(\tilde{u}_r(t,x), J_k)$ denote the total variations of
$u_r(t,x)$, $\tilde{u}_r(t,x)$ on $J_k$ respectively. Suppose
there is a wave $(u_{i-1}(t), u_i(t))$ in $u_r(t,x)$, issued from
$(ks, ir)$ and crosses $J_k$, also $(u_{i-1}, u_i)$ is the
corresponding classical Riemann wave of $(u_{i-1}(t), u_i(t))$
(so $(u_{i-1}, u_i)$ is also issued from $(ks, ir)$ and crosses
$J_k$). If $(u_{i-1}, u_i)$ is a shock wave, then by \eqref{2.7} we can
easily obtain that
\begin{eqnarray} \nonumber
\begin{split}
&|TV((u_{i-1}(t), u_i(t)); J_k) - TV((u_{i-1}, u_i); J_k)|\\
&\leq s
\,  |\frac{\partial{q}}{\partial{u}}( ks, ir, \bar{u}_i)| \, TV((u_{i-1}, u_i);J_k)
     +s \, \big( |q(t_0,x_0,u_{i-1})|+|q(t_0,x_0,u_i)| \big),
\end{split}
\end{eqnarray}
where $\bar{u}_i \in {\cal U}$ and $TV((u_{i-1}(t), u_i(t));
J_k)$, $TV((u_{i-1}, u_i); J_k)$ denote the total variations of
$(u_{i-1}(t), u_i(t))$, $(u_{i-1}, u_i)$ crossing $J_k$.
Similarly, if $(u_{i-1}, u_i) = \bar{u}_i(\xi)$ is a
rarefaction wave with $\xi \in [\xi_1, \xi_2]$, then we obtain
\begin{eqnarray} \nonumber
\begin{split}
&|TV((u_{i-1}(t), u_i(t)); J_k) - TV((u_{i-1}, u_i); J_k)|\\
&\leq s
\, |\frac{\partial{q}}{\partial{u}}( ks, ir, \bar{u}_i(\tilde{\xi})| \,  TV((u_{i-1}, u_i);J_k)
        +s(|q(t_0,x_0,u_{i-1})|+|q(t_0,x_0,u_i)|)
\end{split}
\end{eqnarray}
for some $\tilde{\xi} \in [\xi_1, \xi_2]$ and
$\bar{u}_i(\tilde{\xi}) \in {\cal U}$ . Summing up the previous
inequalities with respect to the waves crossing $J_k$ we obtain
\begin{eqnarray} \nonumber
\begin{split}
& |TV(u_r(J_k)) - TV(\tilde{u}_r(J_k))|
\\
&\leq O(s) \,
\parallel \frac{\partial{q}}{\partial{u}} \parallel_{L^1(\mathbb{R_+} \times
\RR)}
\, \; TV(\tilde{u}_r(J_k))  + O(s) \, \|q\|_{L^1(\mathbb{R_+} \times
\RR)}
\end{split}
\end{eqnarray}
for any mesh curve $J_k$, and this is enough to imply that the
total variations of $u_r(t,x)$ and $\tilde{u}_r(t,x)$ on any mesh
curve $J_k$ are equivalent when $\|q\|_{L^1(\mathbb{R_+} \times
\RR)}$ and $\parallel
\frac{\partial{q}}{\partial{u}} \parallel_{L^1(\mathbb{R_+} \times \RR)}$
are small, we then show the statement.\\

We note that the waves entering each diamond region may come
from two generalized Riemann solutions, we certainly need to know
the constant states of corresponding classical Riemann solutions
at the left and right vertices of diamond region to calculate
those wave strengths separately. We proceed as follows.

First, using the notations in \eqref{4.10}, \eqref{4.11}, we define the
strength of the waves entering the diamond region $\Gamma_{k_0,
h_0},\; k_0 + h_0 \
\rm is \; \rm even$, by
\begin{eqnarray*}
\begin{split}
\eps_* (\Gamma_{k_0, h_0})&:=
|\eps (\tilde{u}_W, u_S ; (k_0 - 1) s, (h_0 - 1) r)| \\
&\quad+|\eps (u_S, \tilde{u}_E ; (k_0 - 1) s, (h_0 + 1) r)|
\end{split}
\label{4.14}
\end{eqnarray*}
\noindent
and the strength of the waves leaving $\Gamma_{k_0, h_0}$ by
\begin{eqnarray}
\eps^* (\Gamma_{k_0, h_0}) :=
|\eps (u_W, \tilde u_N ; k_0 s, h_0 r)| + |\eps
(\tilde{u}_N, u_E ; k_0 s, h_0 r)|.
\label{4.15}
\end{eqnarray}
Since $\tilde{u}_N$ is a constant state in $W_C (u_W, u_E ; k_0
s, h_0 r)$, we can write
\begin{eqnarray}
\eps^* (\Gamma_{k_0, h_0}) =
|\eps (u_W, u_E ; k_0 s, h_0 r)|.
\label{4.16}
\end{eqnarray}
Next,  for $k_0 + h_0 \; \rm is  \;
\rm even$, we let $Q (\Gamma_{k_0, h_0})$ denote the potential of waves interaction in the
diamond $\Gamma_{k_0, h_0}$, i.e.,
\begin{eqnarray*}
Q(\Gamma_{k_0, h_0})
 := D(\eps (\tilde{u}_W,u_S
,(k_0-1)s;(h_0 - 1) r),
\eps(u_S,\tilde{u}_E,(k_0 - 1)s,(h_0+1)r))
\label{4.17}
\end{eqnarray*}
where $D (\cdot,\cdot)$ is defined in \eqref{3.1}. Given a mesh curve
$J$, we note that there are two types of waves crossing $J$. The
first kind of waves are $(\tilde u_{k,h-1},u_{k-1,h})$, $k+h=even$
(waves of type I), the second type of waves are $(u_{k-1,h},\tilde
u_{k,h+1})$, $k+h=even$ (waves of type II). More precisely, waves
of type I are either of the form $(\tilde{u}_{k, h-1}, u_{k-1,
h})$ entering $\Gamma_{k,h}$ (left in-coming waves of
$\Gamma_{k,h}$), or $(\tilde{u}_{k+1, h}, u_{k, h+1})$ leaving
$\Gamma_{k,h}$ (right out-going waves of $\Gamma_{k,h}$). Waves
of type II are either of the form $(u_{k-1, h},
\tilde{u}_{k, h+1})$ entering $\Gamma_{k,h}$ (right in-coming waves of
$\Gamma_{k,h}$), or $(u_{k, h-1},
\tilde{u}_{k+1, h})$ leaving $\Gamma_{k,h}$ (left out-going waves of
$\Gamma_{k,h}$), see Figures 4.2 (a), (b). Next we define the
linear functional $L(J)$ for the waves in $u_r(t,x)$ crossing mesh
curve $J$ by
\begin{eqnarray}
\begin{split}
L (J):= &\sum_{\scriptstyle{\rm type \; I}}
|\eps(\tilde{u}_{k, h - 1}, u_{k - 1, h} ; (k - 1) s, (h - 1) r)|\\
&+\sum_{\scriptstyle{\rm type\;II}} |\eps(u_{k - 1, h},
\tilde{u}_{k, h + 1} ; (k - 1) s, (h + 1) r)|.
\end{split}
\label{4.18}
\end{eqnarray}
\noindent
From previous analysis, we see that functional $L(J)$ is
equivalent to the total variation of $u_r(t,x)$ crossing mesh
curve $J$. Next we define the quadratic functional $Q(J)$ of
$u_r(t,x)$ by
\begin{eqnarray}
Q (J) := \sum_{(\alpha,\beta) } D (\alpha,
\beta)
\label{4.19}
\end{eqnarray}
where the notation $(\alpha,\beta)$ under summation sign denotes
a pair of waves $\alpha$, $\beta$ crossing $J$ and approach, and
$D(\alpha, \beta)$ is given in \eqref{3.1}.
\noindent
Furthermore, we define the $Glimm \; functional \; F(J)$ of
$u_r(t,x)$ for mesh curve $J$ by
\begin{eqnarray}
F (J) := L (J) + K \, Q (J).
\label{4.20}
\end{eqnarray}
\noindent
Our goal is to show that functional $F$ remains uniformly bounded
on all mesh curves provided that constant $K$ in \eqref{4.20} is
sufficiently large, and this leads to the result that functional
$L$ can be bounded by a constant times the total variation of
initial data $u_0(x)$. To show this, the first step is to
estimate the possible changing amount of $L$ and $Q$ when waves
pass through one mesh curve and into an immediate successor. The
estimates of changing amounts of $L$ and $Q$ are stated as
follows.
\begin{center}
\includegraphics[height = 1.5cm, width=10 cm]{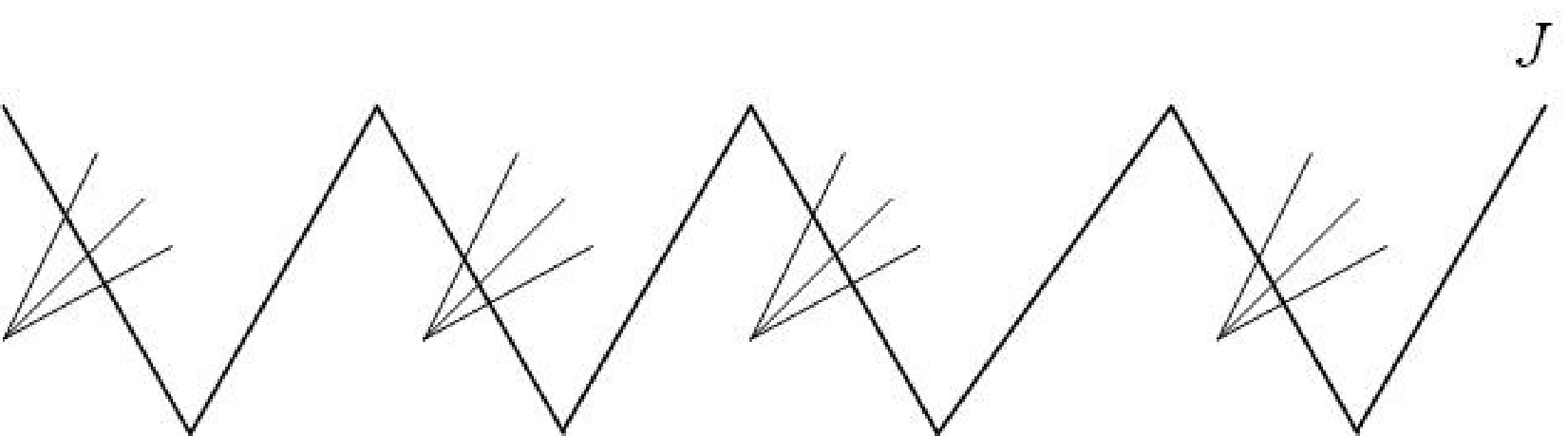} \\
  Figure 4.2(a) :  Waves of type I crossing mesh curve J\\
\includegraphics[height = 1.5cm, width=10 cm]{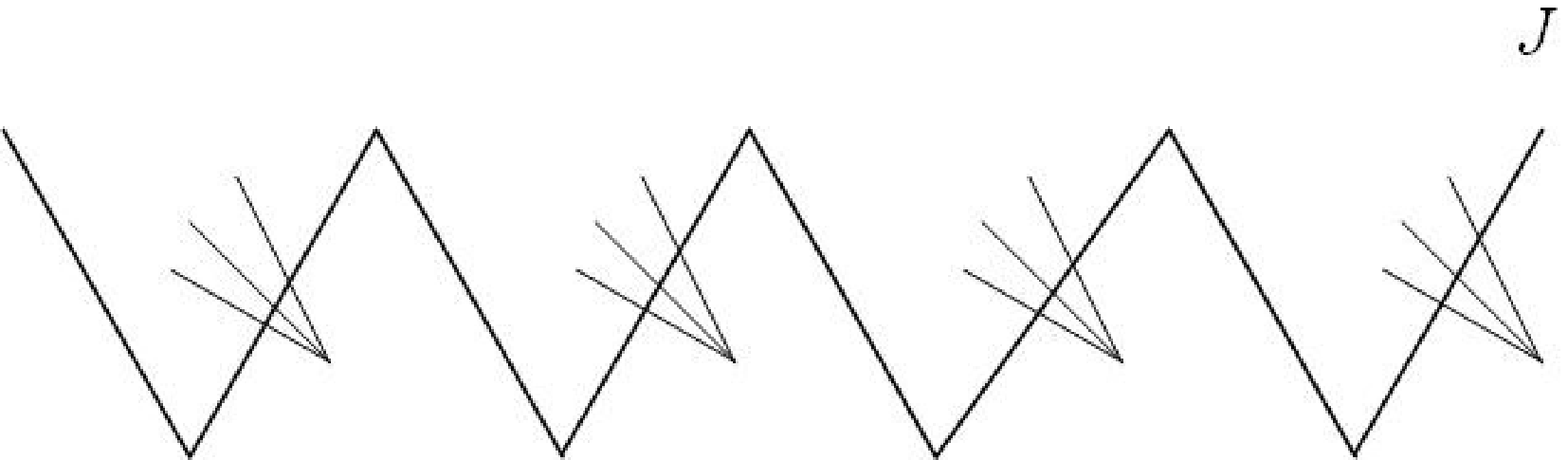} \\
  Figure 4.2(b) :  Waves of type II crossing mesh curve J
\end{center}
\vskip0,5cm
\noindent
{\bf Proposition 4.1.} {\em \quad Given two mesh curves $J_1$ and $J_2$
 such that $J_2$ is an immediate successor of $J_1$, let $\Gamma_{k_0,
h_0}$ denote the diamond region bounded by $J_1$ and $J_2$. Then
functionals $L$ and $Q$ satisfy
\begin{eqnarray}
L (J_2) - L (J_1) = O(1)
\{Q (\Gamma_{k_0, h_0}) + \eps_* (\Gamma_{k_0, h_0})
(C^0_1 + \lambda_* C^0_2 + C^0_3 + C^0_4) s\},
\label{4.21}
\end{eqnarray}
\begin{eqnarray}
\begin{split}
Q (J_2) - Q (J_1)=\;&- Q (\Gamma_{k_0, h_0})
+ O(1) L (J_1)Q(\Gamma_{k_0, h_0}) \\
&+ O(1)  L(J_1) \eps_* (\Gamma_{k_0, h_0}) (C^0_1 +
\lambda_* C^0_2 + C^0_3 + C^0_4)s
\end{split}
\label{4.22}
\end{eqnarray}
 where constants $\lambda_*$ is defined in \eqref{4.1} and  $C^0_j, 1 \leq j \leq 4$, are given by
\begin{eqnarray}
C^0_1 := N^{h_0 r, h_0 r}_{(k_0 - 1) s, k_0 s}(\frac {\partial^2
A}{\partial t \partial u}), \qquad C^0_2 := N^{(h_0 - 1) r, (h_0 +
1)r}_{(k_0 - 1) s, (k_0 - 1) s}(\frac {\partial^2 A }{\partial x
\partial u}),
\label{4.23}
\end{eqnarray}
\begin{eqnarray}
C^0_3 := N^{h_0 r, h_0 r}_{(k_0-1)s, k_0 s} (q(t,x,u)), \qquad
 C^0_4
:= N^{h_0 r, h_0 r}_{(k_0-1)s, k_0 s}(\frac {\partial q}{\partial
u}(t,x,u))
\label{4.24}
\end{eqnarray}
where $N$ is defined in \eqref{3.7}. Note that $\{C^0_j ; 1 \leq j \leq
4\}$ depend on $h_0$,$k_0$.}
\vskip0,5cm
\begin{proof}
\quad Let $u_S, u_W, u_E, u_N$ be the constant states described in
\eqref{4.9}-\eqref{4.10}, we first derive \eqref{4.21}. By the definitions of
$\eps_*$ and $\eps^*$ in \eqref{4.15}, \eqref{4.16} and $L$ in
\eqref{4.18}, we find 
\begin{eqnarray} 
\begin{split}
L (J_2) - L (J_1)&= |\eps (u_W, \tilde{u}_N ; k_0 s, h_0
r)| + |\eps (\tilde{u}_N, u_E ; k_0 s, h_0 r)|\\
&\quad-|\eps (\tilde{u}_W, u_S ; (k_0 - 1) s, (h_0 - 1) r)|\\
&\quad-|\eps (u_S, \tilde{u}_E; (k_0 - 1) s, (h_0 + 1) r)| \\
&=\eps^* (\Delta_{k_0, h_0}) -
\eps_* (\Delta_{k_0, h_0}).
\end{split}
\label{4.25}
\end{eqnarray}
Next, by applying the definition of $\lambda_{*}$ in \eqref{4.1} and the
estimates \eqref{3.32}, \eqref{3.35} to \eqref{4.25} with the choice of $ u_L =
\tilde{u}_W,\ u_M = u_S,\ u_R =
\tilde{u}_E,\
\mu_L = u_W-\tilde{u}_W,\ \mu_R = u_E - \tilde{u}_E$,
 $ t_0 = (k_0 - 1) s$ and $x_0 = h_0 r$, we obtain $$
\eps^* (\Gamma_{k_0, h_0}) =
\eps_* (\Gamma_{k_0, h_0}) + O(1)
Q (\Gamma_{k_0, h_0}) + O(1) \eps_* (\Gamma_{k_0, h_0})
\Big( (C^0_1 + C^0_3 + C^0_4) s + C^0_2 r \Big),
$$ and this gives \eqref{4.21}. \\

To prove \eqref{4.22}, we define several notations for the rest
of the section. First, given $(k,h),\; k+h$= even, we let vector
$\eps_{k-1,h-1/2}$ denote the strength of waves issued from
$((k-1)s,(h-1)r)$ entering $\Gamma_{k,h}$, and let vector
$\eps_{k-1,h+1/2}$ denote the strength of waves issued from
$((k-1)s,(h+1)r)$ entering $\Gamma_{k,h}$. More precisely, the
vector $\eps_{k-1,h-1/2}$ measures the strength of waves of type
I entering $\Gamma_{k,h}$ and $\eps_{k-1,h+1/2}$ measures the
strength of waves of type II entering $\Gamma_{k,h}$. Next, given
a mesh curve $J$, let $J_{[h-1,h]}$ ($J_{[h,h+1]}$ respectively)
denote the segment of $J$ in $\RR_+
\times [(h-1)r,hr]$ ($\RR_+
\times [hr,(h+1)r]$). Then we define vectors $\eps_{J,h-1/2}$, $\eps_{J,h+1/2}$
as the strengths of waves crossing
 $J_{[(h-1),h]}$, $J_{[h,h+1 ]}$ respectively. We will drop the sign $J$ in $\eps_{J,h-1/2}$
 and $\eps_{J,h+1/2}$ when $J$ is specified.
We also set $$
\eps_{W, S} := \eps_{h_0 - 1/2}, \quad
\eps_{S, E} := \eps_{h_0 + 1/2},
$$
 $$\eps_{W, N} := \eps (u_W, \tilde{u}_N ;
k_0 s, h_0 r), \quad \eps_{N, E} :=
\eps (\tilde{u}_N, u_E ; k_0 s, h_0 r).
$$ Since $J_2$ is an immediate successor of $J_1$, the diamond
region bounded by $J_1,\; J_2$ can be specified as $\Gamma(k_0,
h_0)$ with center $(k_0s, \; h_0r)$, and $J_1,\; J_2$ coincide
outside $\Gamma(k_0, h_0)$. We will also drop the signs $J_1,\;
J_2$ without confusion. From the definition of $Q$ in \eqref{4.19}, we
have
\begin{eqnarray*}
\begin{split}
& Q (J_2) - Q (J_1)
\\
& = \sum_{h < h_0 }\Big(D (\eps_{h - 1/2},
\eps_{W, N})
+D (\eps_{h - 1/2}, \eps_{N, E})  -D(\eps_{h - 1/2}, \eps_{W, S}) -D
(\eps_{h -
1/2}, \eps_{S, E})\Big)
 \\
&+ \sum_{h > h_0 + 1} \Big( D (\eps_{W, N}, \eps_{h -
1/2}) +D(\eps_{N, E}, \eps_{h - 1/2})  -D(\eps_{W, S}, \eps_{h - 1/2}) -D
(\eps_{S,
E}, \eps_{h - 1/2})\Big)
\\
&+ D (\eps_{W, N},
\eps_{N, E}) - D (\eps_{W, S}, \eps_{S, E}).
\end{split}
\label{4.26}
\end{eqnarray*}
From \eqref{3.1} we see that
\begin{eqnarray}
D(\eps_{W, N}, \eps_{N, E}) = 0.
\label{4.27}
\end{eqnarray}
Also, for any $h \in \mathbb{Z}$ we observe that
\begin{eqnarray}
D(\eps_{h - 1/2},
\eps_{W, N}) + D (\eps_{h - 1/2}, \eps_{N,
E}) = D (\eps_{h - 1/2}, \eps_{W, E})
\label{4.28}
\end{eqnarray}
for $h < h_0$, and
\begin{eqnarray}
D(\eps_{W, N},
\eps_{h - 1/2}) + D (\eps_{N, E}, \eps_{h -
1/2}) = D (\eps_{W, E},\eps_{h - 1/2})
\label{4.29}
\end{eqnarray}
for $h > h_0 +1$. Thus, by \eqref{4.27}-\eqref{4.29} we obtain
\begin{eqnarray}
\begin{split}
&Q (J_2) - Q (J_1) \\
&=\sum_{h < h_0} \Big( D (\eps_{h - 1/2},\eps_{W,E}) - D (\eps_{h -
1/2},\eps_{W, S})  -D (\eps_{h - 1/2}, \eps_{S, E})\Big)
 \\
& \quad +\sum_{h > h_0 + 1} \Big( D (\eps_{W, E}, \eps_{h - 1/2})-D (\eps_{W, S}, \eps_{h -
1/2}) - D (\eps_{S, E}, \eps_{h - 1/2}) \Big)
 \\
& \quad -D(\eps_{W, S}, \eps_{S, E}).
\end{split}
\label{4.30}
\end{eqnarray}
Finally, applying \eqref{3.33} and \eqref{3.34} to \eqref{4.30} and using the fact
that $D (\eps_{W, S}, \eps_{S, E}) = Q
(\Gamma_{k_0, h_0}),$ we obtain
\begin{eqnarray*}
\begin{split}
&Q (J_2) - Q (J_1) \\
&= -D(\eps_{W, S}, \eps_{S, E}) +\sum_{\scriptstyle h
\in \mathbb{Z} \atop
\scriptstyle{h \neq h_0,h_0+1}} \Big( O(1)|\eps_{h-1/2}| \,
D(\eps_{W,S},\eps_{S,E})\\
& \quad + O(1) |\eps_{h -1/2}| (|\eps_{W, S}|
+|\eps_{S, E}|)  ((C^0_1 + C^0_3 + C^0_4) s + C^0_2 r) \Big)
\\
&=- Q (\Gamma_{k_0, h_0})+ O(1) L(J_1)Q(\Gamma_{k_0, h_0})
+ O(1) L (J_1) \, \eps_* (\Gamma_{k_0, h_0}) ((C^0_1 + C^0_3
+ C^0_4) \, s + C^0_2 \, r), \\
\end{split}
\end{eqnarray*}
which leads to \eqref{4.22}. This completes the proof.
\end{proof}

Before stating a crucial technical lemma, let us introduce a notation about mesh
curves. We say that a mesh curve $J$ is of the type $(k_0, k_0 +
1)$ if all the mesh points on $J$ have the form of \{$(k s, (h +
a_k) r)$ : $k = k_0,\; k_0 + 1$\}.
\vskip0,5cm
\noindent
{\bf Lemma 4.2.} {\em \quad  Given a positive integer $k_0$,
 let $J_1$ and $J_2$ be two mesh curves of type $(k_0 - 1, k_0)$ and
$(k_0, k_0 + 1)$ respectively. We assume that there exists a
positive constant $M_*$ such that
\begin{eqnarray}
L (J_1) \leq M_*.
\label{4.31}
\end{eqnarray}
If $M_*$ is sufficiently small and the constant $K$ in \eqref{4.20} is
sufficiently large, then the functional $F$ satisfies the
following inequality
\begin{eqnarray}
F (J_2) \leq F (J_1) + O(1) s \sum_{h_0 \in \Bbb Z }
\eps_* (\Gamma_{k_0, h_0}) (C^0_1 + \lambda_* C^0_2 +
C^0_3 + C^0_4)
\label{4.32}
\end{eqnarray}
where the bound $O(1)$ depends on $M_*$ and $K$, and the constants
$C^0_j := C^0_j(h_0, k_0)$, $1 \leq j \leq 4$, in \eqref{4.23},
\eqref{4.24} depend on $h_0 \in \mathbb{Z}$.}
\vskip0,5cm
\begin{proof}
\quad  Given $h_0$ $\in$ $\mathbb{Z}$, we multiply
\eqref{4.22} by constant $K$ in \eqref{4.20} and add it to \eqref{4.21}. Then by the
assumption that $J_1$ and $J_2$ are two mesh curves of type $(k_0
- 1, k_0)$ and $(k_0, k_0 + 1)$, we obtain
\begin{eqnarray*}
\begin{split}
F (J_2) - F (J_1) =\;& - K \sum_{h_0 \in \mathbb{Z}} Q
(\Gamma_{k_0, h_0}) + O(1) [1 + K L (J_1)]
\bigg\{ \sum_{h_0 \in \mathbb{Z}} Q (\Gamma_{k_0, h_0})\\
 &+\sum_{h_0 \in \mathbb{Z}} \eps_*
(\Gamma _{k_0, h_0}) (C^0_1 + \lambda_* C^0_2 + C^0_3 + C^0_4) s
\bigg\}.
\end{split}
\end{eqnarray*}
Next, by the observation that $\sum_{h_0 \in \mathbb{Z}} Q
(\Gamma_{k_0, h_0}) = Q (J_1),$ the equation above implies that
\begin{eqnarray*}
\begin{split}
F(J_2) - F(J_1) &= -K Q(J_1)
+ O(1)(1 + K L (J_1))Q(J_1)\\
&\quad+ O(1) s \sum_{h_0 \in \mathbb{Z}} \eps_*
(\Gamma_{k_0, h_0}) (C^0_1 + \lambda_* C^0_2 + C^0_3 + C^0_4)\\
&=Q(J_1)\{K [O(1)L (J_1)-1]+O(1)\}\\
&\quad+ O(1) s \sum_{h_0 \in \mathbb{Z}} \eps_*
(\Gamma_{k_0, h_0}) (C^0_1 + \lambda_* C^0_2 + C^0_3 + C^0_4)\\
&\leq Q(J_1)\{K [O(1)M_*-1]+O(1)\}\\
&\quad+ O(1) s \sum_{h_0 \in \mathbb{Z}} \eps_*
(\Gamma_{k_0, h_0}) (C^0_1 + \lambda_* C^0_2 + C^0_3 + C^0_4).
\end{split}
\label{4.33a}
\end{eqnarray*}
The last inequality is an application of \eqref{4.31}.
We see that the term $K [O(1)M_*-1]+O(1)$ is negative, if $M_*$ is
sufficiently small and $K$ is sufficiently large. Thus, \eqref{4.32}
holds for such $M_*$ and $K$. This completes the proof.
\end{proof}

We now establish the stability of generalized
Glimm method, which is the main result of this section.  We denote by
$TV(\cdot)$ the total variation of a function.

\vskip0,5cm
\noindent
{\bf Theorem 4.3} {\em \quad Fix a constant state $u_*$ and assume
that the initial data $u_0=u_0(x)$ is a function of bounded variation
such that
\begin{eqnarray}
||u_0- u_*||_{L^\infty} \ and \ TV (u_0)\; \hbox{\sl are
sufficiently small}  .
\label{4.34}
\end{eqnarray}
Assume also that the mappings $A(t,x,u) := \frac {Df}{Du} (t, x, u)$
and $q(t,x,u)$ in \eqref{1.6} are smooth and such that 
\begin{eqnarray}
the\; L^1 (\RR_+ \times \RR) \; norm \ of\; \frac
{\partial^2 A}{\partial t
\partial u},\; \lambda_*
\frac{\partial^2 A}{\partial x \partial u}, \ q,\;\frac{\partial
q}{\partial u} \; \hbox{\sl are sufficiently small} .
\label{4.35}
\end{eqnarray}
Then, the approximate solutions $u_r(t,x)$ are bounded uniformly in the $L^\infty$ and $BV$ norms:
\begin{eqnarray}
||u_r - u_*||_{L^\infty (\RR_+ \times \RR)} \leq
O(1) \,
\Big( ||u_0 - u_*||_{L^\infty (\RR)} + TV (u_0) + C\Big),
\label{4.36}
\end{eqnarray}
\begin{eqnarray}
TV (u_r (t, \cdot)) \leq O(1) \, \Big(TV (u_0) + C\Big),
\label{4.37}
\end{eqnarray}
where
\begin{eqnarray}
C := \Big|\Big|\frac {\partial^2 A}{\partial t \partial
u}\Big|\Big|_{L^1 (\RR_+ \times \RR)} + \lambda_* \Big|\Big|\frac {\partial^2
A}{\partial x \partial u}\Big|\Big|_{L^1 (\RR_+ \times \RR)} + ||q||_{L^1 (\RR_+ \times \RR)}
+\Big|\Big|\frac {\partial q}{\partial u}\Big|\Big|_{L^1 (\RR_+ \times \RR)}
\label{4.38}
\end{eqnarray}
Furthermore, the function
$u_r(t,x)$ is Lipschitz continuous in time, i.e., for $t_1,t_2 > 0$,
\begin{eqnarray}
\int_\RR |u_r (t_1, x) - u_r (t_2, x)| \, d x \leq O(1)
(|t_2 - t_1| + s) \, (TV (u_0) + C).
\label{4.39}
\end{eqnarray}
}
\begin{proof}
\quad We apply an induction argument based on Lemma 4.2 to show that the approximate solution $u_r(t,x)$
is uniformly bounded in $L^\infty$ and total variation. First, we
show that the condition \eqref{4.31} in Lemma 4.2 holds under the
assumptions
\eqref{4.34}, \eqref{4.35}. By induction, given $k_0 \in \mathbb{N}$, we let
$J_{k_0 - 1/2}$ denote the mesh curve of type $(k_0 - 1, k_0)$.
For $k_0=1$, we see that
\begin{eqnarray}
F (J_{1/2}) \leq O(1)  \, \big(TV (u_0) + K [TV (u_0)]^2 \big).
\label{4.40}
\end{eqnarray}
This means that there exists a positive constant $M_*$, as
described in \eqref{4.31}, such that $F (J_{1/2}) \leq M_*$, and in
particular, $L(J_{1/2}) \leq M_*$ if $TV(u_0)$ is sufficiently
small. Next, suppose that
\begin{eqnarray}
L (J_{k + 1/2}) \leq M_*  \quad  {\rm for} \ \; k = 0, 1,\ldots,
k_0 - 1.
\label{4.41}
\end{eqnarray}
We intend to show that \eqref{4.41} still holds for $k=k_0$. Since
$J_{k_0 - 1/2}$ is a mesh curve of type $(k_0 - 1, k_0)$, this
implies that $J_{k_0 + 1/2}$ is a mesh curve of type $(k_0, k_0
+1)$ so that Lemma 4.2 can be applied. Therefore we obtain
\begin{eqnarray*}
\begin{split}
F (J_{k_0 + 1/2}) &\leq  F (J_{k_0 - 1/2}) + O(1) s
\sum_{h_0 \in \mathbb{Z}} \eps_* (\Gamma_{k_0, h_0})
(C^0_1 + \lambda_* C^0_2 + C^0_3 + C^0_4)\\
&\ \vdots \\
&\leq  F (J_{1/2}) + O(1) s \sum^{k_0}_{k = 1}
\sum_{h_0
\in
\mathbb{Z}} \eps_* (\Gamma_{k, h_0}) (C^0_1 + \lambda_*
C^0_2 + C^0_3 + C^0_4).
 \end{split}
\end{eqnarray*}
Then, by 
\begin{eqnarray*}
\sum_{h_0 \in \mathbb{Z}} \eps_* (\Gamma_{k, h_0}) =
L (J_{k - 1/2}), \quad  \; k \in \mathbb{N},
\end{eqnarray*}
this leads to
\begin{eqnarray}
F (J_{k_0 + 1/2}) \leq F (J_{1/2}) + O(1)
\sum^{k_0}_{k = 1} s L (J_{k - 1/2}) \sup_{h_0 \in \mathbb{Z}}
(C^0_1 + \lambda_* C^0_2 + C^0_3 + C^0_4).
\label{4.42}
\end{eqnarray}
Next, by \eqref{4.40}-\eqref{4.42} we find 
\begin{eqnarray}
\begin{split}
F(J_{k_0 + 1/2}) \;\leq\; & O(1) (1 + K \, TV (u_0)) \, TV (u_0)\\
 &+O(1) M_* \sum^{k_0 }_{k = 1}  \sup_{h_0 \in \mathbb{Z}}
(C^0_1 + \lambda_* C^0_2 + C^0_3 + C^0_4)s.
\end{split}
\label{4.43}
\end{eqnarray}
From the definitions of $C^0_j$ and the assumption that the
constant $C$ in \eqref{4.38} is finite, we see that
\begin{eqnarray}
\lim_{r \to 0} \sum^{\infty}_{k=1}  \sup_{h_0 \in \mathbb{Z}}
(C^0_1 + \lambda_* C^0_2 + C^0_3 + C^0_4)s = C.
\label{4.44}
\end{eqnarray}
Therefore, from \eqref{4.43} and \eqref{4.44} we obtain the inequality
\begin{eqnarray}
F (J_{k_0 + 1/2}) \leq O(1)
\{(1 + K \, TV (u_0)) \, TV (u_0) + M_* C\},
\label{4.45}
\end{eqnarray}
and in particular,
\begin{eqnarray}
L (J_{k_0 + 1/2}) \leq O(1)
\{(1 + K \, TV (u_0))\, TV (u_0) + M_* C\}.
\label{4.46}
\end{eqnarray}
We note that the functional $L$ in \eqref{4.46} only depends on the
constants $M_*$, $C$ and the total variation of $u_0$, thus it
enables us to choose $TV(u_0)$ and $C$ sufficiently small such
that $O(1) (1 + K \, TV (u_0)) \, TV (u_0) \leq
\frac{M_*}{2}$ and $O(1) \, CM_*\leq
\frac{M_*}{2}$ and this implies that
\begin{eqnarray*}
L (J_{k_0 + 1/2})\leq M_*.
\end{eqnarray*}
Therefore \eqref{4.41} holds for $k = k_0$, we just showed that $L(k_0
+1/2)$ has uniform bound for all $k_0 \in \mathbb{N}$, which
implies that functional $L$ of $u_r(t,x)$ has global bound. Since
$L$ is a functional equivalent to the total variation of
$u_r(t,x)$, we prove that the total variation of $u_r(t,x)$ has
an uniform bound for all $t \geq 0$ and all finite $r>0$, so as
well the $L^\infty$ norm of $u_r(t,x)$.
To prove \eqref{4.37}, we apply \eqref{4.4}, \eqref{4.5} to $u_r(t,x)$ and we use
the fact that $ TV (u_r (k_0 s,\cdot)) = O(1) \, F (J_{k_0 + 1/2})
$ to \eqref{4.45}, then \eqref{4.37} is established. For the proof of \eqref{4.36}
and \eqref{4.39}, we follow the lines of proof in \cite{Gl}. The proof is completed.
\end{proof}
We note that if $\frac {\partial^2 A}{\partial t
\partial u},\; \lambda_*
\frac{\partial^2 A}{\partial x \partial u}, \ q$ and $\frac{\partial
q}{\partial u}$ in \eqref{4.35} belong to $L^\infty$, then inequalities
\eqref{4.36}, \eqref{4.37} and \eqref{4.39} remain valid in a finite interval $[0,
T]$ with $T$ sufficiently small.
\vskip1,0cm
 \setcounter{section}{5}
 \noindent{\bf \large{5. Convergence of the generalized Glimm method}} \\
 \setcounter{equation}{0}

In Section 4 we established the BV stability of the scheme together with a time 
continuity property. 
By Helly's theorem ,there exists a subsequence of approximate solutions, still denoted by $\{u_r(t,x)\}$
and converging strongly in $L^1_{loc}$ to a limit function $u=u(t,x)$. Moreover, by the estimates \eqref{4.31},
\eqref{4.32}, the function $u$ is uniformly bounded and is of bounded variation in $x$. 
We now prove that the limit $u$ is indeed an entropy solution of the Cauchy problem.
The proof relies on the error estimate derived in Section 2. 

\vskip0,5cm
\noindent
{\bf Theorem 5.1}{\em \quad Suppose that the initial data $u_0(x)$
is sufficiently close to a constant state in $L^\infty$ and $BV$,
and that the $L^1$ norms of $\frac {\partial^2 A}{\partial x
\partial u}$,\;$\frac{\partial^2 A}{\partial t
\partial x}$,\;$q$, and $\frac {\partial q}{\partial
u}$ are sufficiently small in $\RR_{+}
\times \RR$. Let $\{u_r(t,x): r >0 \}
$ be the sequence of approximate solutions constructed by the
generalized Glimm scheme (4.3)-(4.5). Then, for any
equidistributed sequence $\{a_k\}_{k \in
\mathbb{N}}$, there exists a subsequence of $\{u_r(t,x)\}
$ converging in $L^1_{loc}$ to a function $u=u(t,x)$ which is an
entropy solution of the Cauchy problem (1.1), (1.2).}
\vskip0,5cm

\

\noindent{\bf Remark 5.2.}  Assume that ${\cal U}$ is a
convex subset of $R^p$, we say that $(U,F)$, $U : {\cal U} \in
\RR^p
\to
\RR$ and $F : \RR_+ \times
\RR \times {\cal U} \to \RR$, is an entropy pair of the system \eqref{1.1}
if $U $ is a convex function on ${\cal U}$ and
\begin{eqnarray*}
\frac {\partial F}{\partial u} = \frac{D U}{D u}\; \frac{\partial f}{\partial u} \;\;\;\;\; {\rm on}
\;
\RR_+ \times \RR \times {\cal U}.
\label{5.1}
\end{eqnarray*}
Furthermore, a function $u : \RR_+ \times \RR  \to
\RR^p$ is called an entropy solution of \eqref{1.1} if $u=u(t,x)$ is a weak solution
of \eqref{1.1} satisfying
\begin{eqnarray}
\partial_t U (u) + \partial_x (F (t, x, u)) \leq \frac {D U}{D u} (u) \{g (t, x,
u)-(\partial_xf)(t,x,u)\} +(\partial_x F)(t,x,u)
\label{5.2}
\end{eqnarray}
in the sense of distributions, for every entropy pair $(U, F)$.

\

\begin{proof} \quad The proof is based on the result of Proposition~2.1. 
Let $\{u_r(t,x)\}$ denote a sequence of approximate solutions
constructed by generalized Glimm scheme \eqref{4.3}-\eqref{4.5}.
Then, by the stability result and Helly's theorem, there exists a
subsequence of $\{u_r(t,x)\}$ converging almost everywhere to a
function $u \in L^1_{loc}$ with bounded total variation.
Given any test-function $\theta :
\RR_+ \times
\RR \to \RR$ with we define the residual of $u_r$ as
\begin{eqnarray*}
{\cal R}(u_r, \theta) :=  \int_{\RR_+}
\int_\RR
\{u_r \partial_t \theta + f (t, x, u_r) \partial_x \theta +
g (t, x, u_r) \theta\} d x \ d t + \int^{\infty}_{-\infty} u_0(x)
\theta(0,x) d x
\label{5.3}
\end{eqnarray*}
Note that $u$ is a weak solution to \eqref{1.1},
\eqref{1.2} if and only if ${\cal R}(u, \theta)=0$ for any test-function
$\theta$. By Lebesgue's theorem, we see that
\begin{eqnarray*}
|{\cal R}(u_r, \theta) - {\cal R}(u,\theta)| \rightarrow 0
\quad
{\rm as} \ r
\to 0.
\label{5.4}
\end{eqnarray*}
Thus, to show that $u$ is a weak solution of \eqref{1.1},
\eqref{1.2}, it is equivalent to show that ${\cal R}(u_r(t,x),
\theta)$ tends to zero as $r$ vanishes. To show this, we
first let $\chi^{k_0, h_0}_{supp \;\; (\theta)}$ denote a
characteristic function having the same support as the test-function $\theta$. Then, by 
construction of $u_r(t,x)$ and by 
\eqref{2.11}, we can write
\begin{eqnarray}
\begin{split}
&{\cal R} (u_r, \theta)\\
=& \sum_{k_0=0}^{\infty} \sum_{h_0+k_0  \atop \rm even}
\int^{(k_0 + 1) s}_{k_0 s}\int^{(h_0 + 1) r}_{(h_0 - 1) r}
\Big( u_r \partial_t \theta + f (t, x, u_r) \partial x \theta + g (t, x, u_r) \theta\Big) \, d x \ d t\\
=& \sum_{k_0=0}^{\infty} \sum_{h_0+k_0  \atop \rm even} O(1) (s^2
+ r^2) (s + r + |u_{k_0, h_0 + 1} - u_{k_0, h_0 - 1}|)
\chi^{k_0, h_0}_{supp \;\; (\theta)}\\
&+\Bigg( \sum_{k_0=0}^{\infty} \sum_{h_0+k_0  \atop \rm even}
  \bigg( \int^{(h_0 + 1) r}_{(h_0 - 1) r} u_r ((k_0 + 1)
s -, x) \theta ((k_0 + 1) s, x) d x\\
& \quad - \int^{(h_0 + 1) r}_{(h_0 - 1) r} u_r (k_0 s +, x)\theta (k_0
s, x) d x \bigg) +\int^{\infty}_{-\infty} u_0(x) \theta(0,x) d x
\Bigg) \\ &+ \sum_{k_0=0}^{\infty}
\sum_{h_0+k_0
\atop \rm even}
\Bigg( \int^{(k_0 + 1) s}_{k_0 s} f (t, (h_0 + 1) r,
u_r (t, (h_0 + 1) r - )) \theta (t, (h_0 + 1) r)  d t \\
& \quad -\int^{(k_0 + 1) s}_{k_0 s} f (t, (h_0 - 1) r, u_r (t, (h_0 - 1) r
+ )) \theta (t, (h_0 - 1) r)  d t
\Bigg).
\end{split}
\label{5.5}
\end{eqnarray}
\noindent
Let $\Omega_1(r), \; \Omega_2(r)$ and $\Omega_3(r)$ denote the
terms on the right hand side of \eqref{5.5}, respectively. We
first estimate $\Omega_1(r)$. By a direct calculation and
\eqref{4.2}, \eqref{4.37}, we obtain
\begin{eqnarray}
\begin{split}
\Omega_1(r) =& O(1) r +\sum_{k_0=0}^{\infty} \sum_{h_0+k_0  \atop \rm even} O(1) (s^2 + r^2)
( |u_{k_0, h_0 + 1} - u_{k_0, h_0 - 1}|)
\chi^{k_0, h_0}_{supp \;\; (\theta)}\\
\leq&O(1)r+\sum_{k_0 \in \mathbb{N}}O(1) (s^2 +
r^2)(T.V.\{u_0(x)\}+C)\chi^{k_0, h_0}_{supp \;\; (\theta)}\\
\leq&O(1)r.
\end{split}
\label{5.6}
\end{eqnarray}
Next we calculate $\Omega_3(r)$. By the property of the Lipschitz
continuity of $f$, $q$ and \eqref{2.7},
\eqref{4.5}, we obtain
\begin{eqnarray*}
\begin{split}
&\Omega_3(r)
\\
 &=O(1) \sum_{k_0=0}^{\infty} \sum_{h_0+k_0  \atop \rm even}
\int^{(k_0 + 1) s}_{k_0 s} |u_r (t, (h_0 + 1) r + ) - u_r (t,
(h_0 - 1) r -)|\cdot (\chi^{k_0, h_0}_{supp \;\; (\theta)}) d t \\
&=O(1) \sum_{k_0=0}^{\infty} \sum_{h_0+k_0  \atop \rm even}\Big(
\int^{(k_0 + 1) s}_{k_0 s} t |q (k_0 s, (h_0 + 2) r,
\tilde{u}_{k_0, h_0 + 1}) - q (k_0 s, h_0 r,
\tilde{u}_{k_0, h_0
+ 1})|\cdot (\chi^{k_0, h_0}_{supp \;\; (\theta)}) d t \Big)
\\
&=O(1) \sum_{k_0}  \sum_{h_0}
 \int^{(k_0 + 1) s}_{k_0 s} t \  r \  \cdot (\chi^{k_0, h_0}_{supp \;\;
 (\theta)}) \ d t.
\end{split}
\end{eqnarray*}
It follows that
\begin{eqnarray}
\Omega_3(r) = O(1) r.
\label{5.7}
\end{eqnarray}
It remains to estimate $\Omega_2(r)$. It is a standard matter to check that
\begin{eqnarray*}  
\begin{split}
\Omega_2(r) =&-\sum_{k_0=1}^{\infty} \sum_{h_0+k_0  \atop \rm even}\int^{(h_0+1)r}_{(h_0-1)r}[u_r] (k_0 s, x)
\theta (k_0 s, x) d x\\
&-\int^{+ \infty}_{-
\infty}(u_r(0,x)-u_0(x))\theta(0,x)d x
\end{split}
\label{5.8}
\end{eqnarray*}
where $ [u_r](k_0 s,x):= u_r(k_0 s+,x)-u_r(k_0 s-,x)$. We let
$J(\{a_k\},r,\theta)$ denote the term
$$
\displaystyle
\sum_{k_0=1}^{\infty} \sum_{h_0+k_0  \atop \rm even}\int^{(h_0+1)r}_{(h_0-1)r}[u_r] (k_0 s, x)
\theta (k_0 s, x) d x.
$$
By the construction of $u_r(t,x)$ in \eqref{4.3}, we see that the term $\int^{+ \infty}_{-
\infty}(u_r(0,x)-u_0(x))\theta(0,x)d x$ on the right hand side of \eqref{5.8} vanishes as $r$ tends to
zero. In addition, by a result of Liu \cite{Li1} we obtain that, for
any equidistributed sequence $\{a_k\}_{k \in
\mathbb{N}}$,
$J(\{a_k\},r,\theta)$ tends to zero as $r$ approaches to zero.
This implies that
\begin{eqnarray}
\Omega_2(r) \to 0 \quad {\rm as} \ r \to 0
\label{5.9}
\end{eqnarray}
for every equidistributed sequence $\{a_k\}_{k \in
\mathbb{N}}$. We refer the reader to \cite{Li1} for the details of the estimate of $\Omega_2(r)$. Finally, by \eqref{5.6},
\eqref{5.7} and
\eqref{5.9}, we obtain
\begin{eqnarray*}
{\cal R}(u_r, \theta) \to 0 \; {\rm in} \; L^1 \quad {\rm
as} \ r
\to 0,
\end{eqnarray*}
which means that the limit function $u$ satisfies $ {\cal
R}(u,
\theta) = 0$. Therefore, $u$ is a weak solution of the Cauchy problem \eqref{1.1}, \eqref{1.2}.

To prove that $u$ is an entropy solution satisfying the entropy
inequality \eqref{5.2}, it is equivalent to show that, for any entropy
pair $(U,F)$ and test-function $\theta \geq 0$, the function $u$ satisfies
\begin{eqnarray}
\int_{\RR_+}\int_\RR U(u)\theta_t+F(t,x,u)\theta_x+P(t,x,u)\theta d x d
t+\int_\RR U(u_0(x)) \, \theta(x,0) \, d x \geq 0,
\label{5.10}
\end{eqnarray}
with
\begin{eqnarray*}
P(t,x,u) := \frac{D U}{D u}\cdot (g-\frac{\partial f}{\partial
x})(t,x,u)+(\partial_x F)(t,x,u).
\end{eqnarray*}
We note that the result of Proposition~2.1 can be applied to show
that $u(t,x)$ satisfies \eqref{5.10} for any entropy pair $(U,F)$. In turn,
this implies that $u$ is an entropy solution of the Cauchy problem
\eqref{1.1}, \eqref{1.2}, and the proof of Theorem 5.1 is completed.
 \\
\end{proof}


 \noindent{\bf \large{Acknowledgments.}} \quad
This research was partially supported by the National Science Council (NSC) of Taiwan,
and by the Centre National de la Recherche Scientifique (CNRS).


\end{document}